\documentclass[11pt, a4paper]{article}
\usepackage{a4wide}
\usepackage{amsmath, amsthm}
\usepackage[all,cmtip]{xy}
\entrymodifiers={+!!<0pt,\fontdimen22\textfont2>}

\usepackage{color}
\usepackage{ifpdf}
\ifpdf 
  \usepackage[pdftex]{graphicx}
  \DeclareGraphicsExtensions{.pdf,.png,.mps}
  \usepackage{pgf}
  \usepackage{tikz}
\else 
  \usepackage{graphicx}
  \DeclareGraphicsExtensions{.eps,.bmp}
  \usepackage{pgf}
  \usepackage{tikz}
  \usepackage{pstricks}
\fi
\usepackage{epic}
\usepackage{wrapfig}

\newtheorem{thm}{Theorem}
\newtheorem{cor}[thm]{Corollary}
\newtheorem{lem}[thm]{Lemma}

\theoremstyle{definition}

\newtheorem*{rmk}{Remark}
\newtheorem*{ex}{Example}

\newcommand{\nt}{\par\noindent{\it Note. }}
\newcommand{\qu}{\par\noindent{\it Question. }}
\newcommand{\subtitle}[1]{\par\medskip\noindent{\bf #1.}}

\newcommand{\abs}[1]{\left\vert#1\right\vert}
\newcommand{\set}[1]{\left\{#1\right\}}
\newcommand\qbin[2]{{#1 \brack #2}_{q}}
\newcommand\qmultinomial[2]{{#1 \brack #2}_{q}}
\def\area{{\mathop {\rm area}}}

\title{An involution for symmetry of hook length and part length of pointed partitions}
\author{Heesung Shin and Jiang Zeng}

\begin{document}
\maketitle
\begin{abstract}
A {\em pointed partition} of $n$ is a pair $(\lambda, v)$ where $\lambda\vdash n$ and $v$ is a cell in its Ferrers diagram.
We construct an involution on pointed partitions of $n$ exchanging ``hook length'' and ``part length''. This gives a bijective proof of a recent result of Bessenrodt and Han.
\end{abstract}

\section{Introduction}
A partition $\lambda$ is a sequence of positive integers $\lambda = (\lambda_1, \lambda_2, \ldots, \lambda_{\ell})$ such that $\lambda_1 \geq \lambda_2 \geq \cdots \geq \lambda_\ell > 0$. The integers $\lambda_1, \ldots, \lambda_{\ell}$ are called the parts of $\lambda$, the number $\ell$ of parts is denoted by $\ell(\lambda)$ and called the length of $\lambda$. The sum of parts is denoted by $\abs{\lambda}$. If $\abs{\lambda}=n$, we say that $\lambda$ is a partition of $n$ and write $\lambda \vdash n$.
Each partition can be represented by its {\em Ferrers diagram}, we shall identify a partition with its Ferrers diagram. We refer the reader to Andrews' book \cite{MR1634067} for general reference on partitions.

\begin{figure}[b]
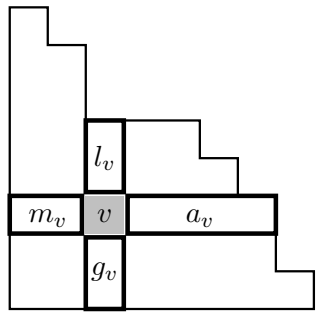

$$
\centering
\begin{pgfpicture}{58.00mm}{48.00mm}{102.00mm}{92.00mm}
\pgfsetxvec{\pgfpoint{1.00mm}{0mm}}
\pgfsetyvec{\pgfpoint{0mm}{1.00mm}}
\color[rgb]{0,0,0}\pgfsetlinewidth{0.30mm}\pgfsetdash{}{0mm}
\color[rgb]{0.75294,0.75294,0.75294}\pgfmoveto{\pgfxy(70.00,60.00)}\pgflineto{\pgfxy(75.00,60.00)}\pgflineto{\pgfxy(75.00,65.00)}\pgflineto{\pgfxy(70.00,65.00)}\pgfclosepath\pgffill
\color[rgb]{0,0,0}\pgfputat{\pgfxy(72.50,61.50)}{\pgfbox[bottom,left]{\fontsize{11.38}{13.66}\selectfont \makebox[0pt]{$v$}}}
\pgfputat{\pgfxy(65.00,61.50)}{\pgfbox[bottom,left]{\fontsize{11.38}{13.66}\selectfont \makebox[0pt]{$m_v$}}}
\pgfputat{\pgfxy(85.00,61.50)}{\pgfbox[bottom,left]{\fontsize{11.38}{13.66}\selectfont \makebox[0pt]{$a_v$}}}
\pgfputat{\pgfxy(72.50,69.00)}{\pgfbox[bottom,left]{\fontsize{11.38}{13.66}\selectfont \makebox[0pt]{$l_v$}}}
\pgfputat{\pgfxy(72.50,54.50)}{\pgfbox[bottom,left]{\fontsize{11.38}{13.66}\selectfont \makebox[0pt]{$g_v$}}}
\pgfmoveto{\pgfxy(60.00,50.00)}\pgflineto{\pgfxy(100.00,50.00)}\pgflineto{\pgfxy(100.00,55.00)}\pgflineto{\pgfxy(95.00,55.00)}\pgflineto{\pgfxy(95.00,65.00)}\pgflineto{\pgfxy(90.00,65.00)}\pgflineto{\pgfxy(90.00,70.00)}\pgflineto{\pgfxy(85.00,70.00)}\pgflineto{\pgfxy(85.00,75.00)}\pgflineto{\pgfxy(70.00,75.00)}\pgflineto{\pgfxy(70.00,85.00)}\pgflineto{\pgfxy(65.00,85.00)}\pgflineto{\pgfxy(65.00,90.00)}\pgflineto{\pgfxy(60.00,90.00)}\pgfclosepath\pgfstroke
\pgfsetlinewidth{0.60mm}\pgfmoveto{\pgfxy(60.00,60.00)}\pgflineto{\pgfxy(69.50,60.00)}\pgflineto{\pgfxy(69.50,65.00)}\pgflineto{\pgfxy(60.00,65.00)}\pgfclosepath\pgfstroke
\pgfmoveto{\pgfxy(75.50,60.00)}\pgflineto{\pgfxy(95.00,60.00)}\pgflineto{\pgfxy(95.00,65.00)}\pgflineto{\pgfxy(75.50,65.00)}\pgfclosepath\pgfstroke
\pgfmoveto{\pgfxy(70.00,50.00)}\pgflineto{\pgfxy(75.00,50.00)}\pgflineto{\pgfxy(75.00,59.50)}\pgflineto{\pgfxy(70.00,59.50)}\pgfclosepath\pgfstroke
\pgfmoveto{\pgfxy(70.00,65.50)}\pgflineto{\pgfxy(75.00,65.50)}\pgflineto{\pgfxy(75.00,75.00)}\pgflineto{\pgfxy(70.00,75.00)}\pgfclosepath\pgfstroke
\end{pgfpicture}%
$$
\caption{The arm, leg, coarm and  coleg length of $(\lambda, v)$}
\label{fig:pointed}
\end{figure}

\begin{figure}[t]
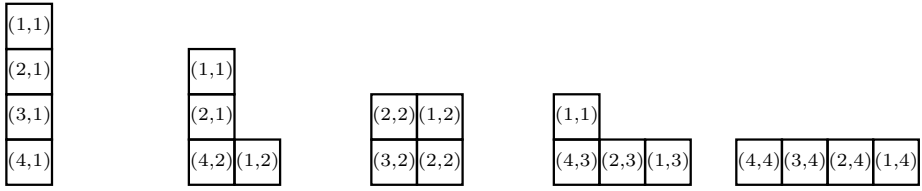

$$
\begin{pgfpicture}{4.00mm}{33.85mm}{128.00mm}{61.85mm}
\pgfsetxvec{\pgfpoint{0.60mm}{0mm}}
\pgfsetyvec{\pgfpoint{0mm}{0.60mm}}
\color[rgb]{0,0,0}\pgfsetlinewidth{0.30mm}\pgfsetdash{}{0mm}
\pgfmoveto{\pgfxy(10.00,89.74)}\pgflineto{\pgfxy(20.00,89.74)}\pgflineto{\pgfxy(20.00,99.74)}\pgflineto{\pgfxy(10.00,99.74)}\pgfclosepath\pgfstroke
\pgfmoveto{\pgfxy(50.00,79.74)}\pgflineto{\pgfxy(60.00,79.74)}\pgflineto{\pgfxy(60.00,89.74)}\pgflineto{\pgfxy(50.00,89.74)}\pgfclosepath\pgfstroke
\pgfmoveto{\pgfxy(10.00,79.74)}\pgflineto{\pgfxy(20.00,79.74)}\pgflineto{\pgfxy(20.00,89.74)}\pgflineto{\pgfxy(10.00,89.74)}\pgfclosepath\pgfstroke
\pgfmoveto{\pgfxy(10.00,69.74)}\pgflineto{\pgfxy(20.00,69.74)}\pgflineto{\pgfxy(20.00,79.74)}\pgflineto{\pgfxy(10.00,79.74)}\pgfclosepath\pgfstroke
\pgfmoveto{\pgfxy(10.00,59.74)}\pgflineto{\pgfxy(20.00,59.74)}\pgflineto{\pgfxy(20.00,69.74)}\pgflineto{\pgfxy(10.00,69.74)}\pgfclosepath\pgfstroke
\pgfmoveto{\pgfxy(50.00,69.74)}\pgflineto{\pgfxy(60.00,69.74)}\pgflineto{\pgfxy(60.00,79.74)}\pgflineto{\pgfxy(50.00,79.74)}\pgfclosepath\pgfstroke
\pgfmoveto{\pgfxy(50.00,59.74)}\pgflineto{\pgfxy(60.00,59.74)}\pgflineto{\pgfxy(60.00,69.74)}\pgflineto{\pgfxy(50.00,69.74)}\pgfclosepath\pgfstroke
\pgfmoveto{\pgfxy(60.00,59.74)}\pgflineto{\pgfxy(70.00,59.74)}\pgflineto{\pgfxy(70.00,69.74)}\pgflineto{\pgfxy(60.00,69.74)}\pgfclosepath\pgfstroke
\pgfmoveto{\pgfxy(90.00,69.74)}\pgflineto{\pgfxy(100.00,69.74)}\pgflineto{\pgfxy(100.00,79.74)}\pgflineto{\pgfxy(90.00,79.74)}\pgfclosepath\pgfstroke
\pgfmoveto{\pgfxy(90.00,59.74)}\pgflineto{\pgfxy(100.00,59.74)}\pgflineto{\pgfxy(100.00,69.74)}\pgflineto{\pgfxy(90.00,69.74)}\pgfclosepath\pgfstroke
\pgfmoveto{\pgfxy(100.00,69.74)}\pgflineto{\pgfxy(110.00,69.74)}\pgflineto{\pgfxy(110.00,79.74)}\pgflineto{\pgfxy(100.00,79.74)}\pgfclosepath\pgfstroke
\pgfmoveto{\pgfxy(100.00,59.74)}\pgflineto{\pgfxy(110.00,59.74)}\pgflineto{\pgfxy(110.00,69.74)}\pgflineto{\pgfxy(100.00,69.74)}\pgfclosepath\pgfstroke
\pgfmoveto{\pgfxy(130.00,69.74)}\pgflineto{\pgfxy(140.00,69.74)}\pgflineto{\pgfxy(140.00,79.74)}\pgflineto{\pgfxy(130.00,79.74)}\pgfclosepath\pgfstroke
\pgfmoveto{\pgfxy(130.00,59.74)}\pgflineto{\pgfxy(140.00,59.74)}\pgflineto{\pgfxy(140.00,69.74)}\pgflineto{\pgfxy(130.00,69.74)}\pgfclosepath\pgfstroke
\pgfmoveto{\pgfxy(140.00,59.74)}\pgflineto{\pgfxy(150.00,59.74)}\pgflineto{\pgfxy(150.00,69.74)}\pgflineto{\pgfxy(140.00,69.74)}\pgfclosepath\pgfstroke
\pgfmoveto{\pgfxy(150.00,59.74)}\pgflineto{\pgfxy(160.00,59.74)}\pgflineto{\pgfxy(160.00,69.74)}\pgflineto{\pgfxy(150.00,69.74)}\pgfclosepath\pgfstroke
\pgfmoveto{\pgfxy(170.00,59.74)}\pgflineto{\pgfxy(180.00,59.74)}\pgflineto{\pgfxy(180.00,69.74)}\pgflineto{\pgfxy(170.00,69.74)}\pgfclosepath\pgfstroke
\pgfmoveto{\pgfxy(180.00,59.74)}\pgflineto{\pgfxy(190.00,59.74)}\pgflineto{\pgfxy(190.00,69.74)}\pgflineto{\pgfxy(180.00,69.74)}\pgfclosepath\pgfstroke
\pgfmoveto{\pgfxy(190.00,59.74)}\pgflineto{\pgfxy(200.00,59.74)}\pgflineto{\pgfxy(200.00,69.74)}\pgflineto{\pgfxy(190.00,69.74)}\pgfclosepath\pgfstroke
\pgfmoveto{\pgfxy(200.00,59.74)}\pgflineto{\pgfxy(210.00,59.74)}\pgflineto{\pgfxy(210.00,69.74)}\pgflineto{\pgfxy(200.00,69.74)}\pgfclosepath\pgfstroke
\pgfputat{\pgfxy(15.00,94.00)}{\pgfbox[bottom,left]{\fontsize{6.83}{8.19}\selectfont \makebox[0pt]{(1,1)}}}
\pgfputat{\pgfxy(55.00,84.00)}{\pgfbox[bottom,left]{\fontsize{6.83}{8.19}\selectfont \makebox[0pt]{(1,1)}}}
\pgfputat{\pgfxy(55.00,74.00)}{\pgfbox[bottom,left]{\fontsize{6.83}{8.19}\selectfont \makebox[0pt]{(2,1)}}}
\pgfputat{\pgfxy(55.00,64.00)}{\pgfbox[bottom,left]{\fontsize{6.83}{8.19}\selectfont \makebox[0pt]{(4,2)}}}
\pgfputat{\pgfxy(15.00,84.00)}{\pgfbox[bottom,left]{\fontsize{6.83}{8.19}\selectfont \makebox[0pt]{(2,1)}}}
\pgfputat{\pgfxy(15.00,74.00)}{\pgfbox[bottom,left]{\fontsize{6.83}{8.19}\selectfont \makebox[0pt]{(3,1)}}}
\pgfputat{\pgfxy(15.00,64.00)}{\pgfbox[bottom,left]{\fontsize{6.83}{8.19}\selectfont \makebox[0pt]{(4,1)}}}
\pgfputat{\pgfxy(65.00,64.00)}{\pgfbox[bottom,left]{\fontsize{6.83}{8.19}\selectfont \makebox[0pt]{(1,2)}}}
\pgfputat{\pgfxy(95.00,74.00)}{\pgfbox[bottom,left]{\fontsize{6.83}{8.19}\selectfont \makebox[0pt]{(2,2)}}}
\pgfputat{\pgfxy(105.00,74.00)}{\pgfbox[bottom,left]{\fontsize{6.83}{8.19}\selectfont \makebox[0pt]{(1,2)}}}
\pgfputat{\pgfxy(95.00,64.00)}{\pgfbox[bottom,left]{\fontsize{6.83}{8.19}\selectfont \makebox[0pt]{(3,2)}}}
\pgfputat{\pgfxy(105.00,64.00)}{\pgfbox[bottom,left]{\fontsize{6.83}{8.19}\selectfont \makebox[0pt]{(2,2)}}}
\pgfputat{\pgfxy(135.00,74.00)}{\pgfbox[bottom,left]{\fontsize{6.83}{8.19}\selectfont \makebox[0pt]{(1,1)}}}
\pgfputat{\pgfxy(155.00,64.00)}{\pgfbox[bottom,left]{\fontsize{6.83}{8.19}\selectfont \makebox[0pt]{(1,3)}}}
\pgfputat{\pgfxy(145.00,64.00)}{\pgfbox[bottom,left]{\fontsize{6.83}{8.19}\selectfont \makebox[0pt]{(2,3)}}}
\pgfputat{\pgfxy(135.00,64.00)}{\pgfbox[bottom,left]{\fontsize{6.83}{8.19}\selectfont \makebox[0pt]{(4,3)}}}
\pgfputat{\pgfxy(175.00,64.00)}{\pgfbox[bottom,left]{\fontsize{6.83}{8.19}\selectfont \makebox[0pt]{(4,4)}}}
\pgfputat{\pgfxy(185.00,64.00)}{\pgfbox[bottom,left]{\fontsize{6.83}{8.19}\selectfont \makebox[0pt]{(3,4)}}}
\pgfputat{\pgfxy(195.00,64.00)}{\pgfbox[bottom,left]{\fontsize{6.83}{8.19}\selectfont \makebox[0pt]{(2,4)}}}
\pgfputat{\pgfxy(205.00,64.00)}{\pgfbox[bottom,left]{\fontsize{6.83}{8.19}\selectfont \makebox[0pt]{(1,4)}}}
\end{pgfpicture}%
$$
\caption{The distribution of $(h_v,p_v)$ on $\mathcal{F}_4$}
\label{fig:part}
\end{figure}

A {\em pointed partition} of $n$ is a pair $(\lambda, v)$ where $\lambda\vdash n$ and $v$ is a cell in its Ferrers diagram. Let $\mathcal{F}_n$ be the set of pointed partitions of $n$.
For each pointed partition $(\lambda, v)$, we define the arm length $a_v:=a_v(\lambda)$ (resp. leg length $l_v$, coarm length $m_v$, coleg length $g_v$) to be the number of cells  lying  in the same row as $v$ and to the right of $v$ (resp. in the same column as $v$ and above $v$, in the same row as $v$ and to the left of $v$, in the same column as $v$ and under $v$), see Figure~\ref{fig:pointed}. The \emph{hook length}  $h_v $ and \emph{part length} $p_v$ of $v$ in $\lambda$ are defined by $h_v = l_v + a_v + 1$ and $p_v=m_v + a_v+1$, respectively. The joint distribution of $(h_v, p_v)$ on ${\cal F}_4$ is given in Figure~\ref{fig:part}, where $(h_v,p_v)$ is written in $(\lambda, v)$.

In a recent paper, generalizing the result of Bessenrodt \cite{MR1682922}, Bessenrodt and Han \cite{BESSENRODT:2009:HAL-00395700:1}  showed that the bivariate joint distribution $(h_v, p_v)$ over $\mathcal{F}_n$ is {\em symmetric}, that is,
$$
\sum_{(\lambda,v)\in \mathcal{F}_n} x^{h_v} y^{p_v} = \sum_{(\lambda,v)\in \mathcal{F}_n} x^{p_v} y^{h_v}.
$$
Since the proof in \cite{BESSENRODT:2009:HAL-00395700:1} uses a generating function argument, this raises the natural question of finding a bijective proof of their result.
The aim of this paper is to construct an involution $\Phi$ on $\mathcal{F}_n$ exchanging hook length and part length and give bijective proofs of the main results in \cite{BESSENRODT:2009:HAL-00395700:1}.

\section{Hook and rim hook}
\begin{figure}
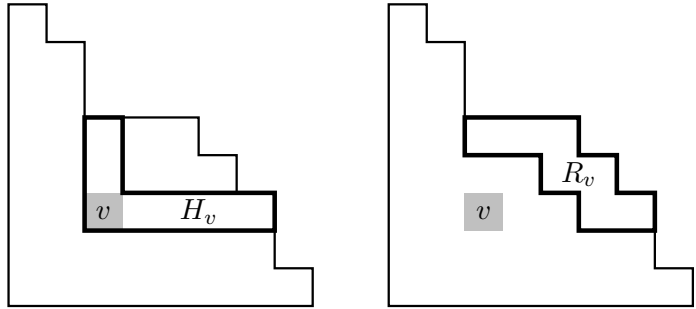

$$
\centering
\begin{pgfpicture}{108.00mm}{48.00mm}{202.00mm}{92.00mm}
\pgfsetxvec{\pgfpoint{1.00mm}{0mm}}
\pgfsetyvec{\pgfpoint{0mm}{1.00mm}}
\color[rgb]{0,0,0}\pgfsetlinewidth{0.30mm}\pgfsetdash{}{0mm}
\color[rgb]{0.75294,0.75294,0.75294}\pgfmoveto{\pgfxy(120.00,60.00)}\pgflineto{\pgfxy(125.00,60.00)}\pgflineto{\pgfxy(125.00,65.00)}\pgflineto{\pgfxy(120.00,65.00)}\pgfclosepath\pgffill
\color[rgb]{0,0,0}\pgfputat{\pgfxy(122.50,61.50)}{\pgfbox[bottom,left]{\fontsize{11.38}{13.66}\selectfont \makebox[0pt]{$v$}}}
\pgfmoveto{\pgfxy(110.00,50.00)}\pgflineto{\pgfxy(150.00,50.00)}\pgflineto{\pgfxy(150.00,55.00)}\pgflineto{\pgfxy(145.00,55.00)}\pgflineto{\pgfxy(145.00,65.00)}\pgflineto{\pgfxy(140.00,65.00)}\pgflineto{\pgfxy(140.00,70.00)}\pgflineto{\pgfxy(135.00,70.00)}\pgflineto{\pgfxy(135.00,75.00)}\pgflineto{\pgfxy(120.00,75.00)}\pgflineto{\pgfxy(120.00,85.00)}\pgflineto{\pgfxy(115.00,85.00)}\pgflineto{\pgfxy(115.00,90.00)}\pgflineto{\pgfxy(110.00,90.00)}\pgfclosepath\pgfstroke
\color[rgb]{0.75294,0.75294,0.75294}\pgfmoveto{\pgfxy(170.00,60.00)}\pgflineto{\pgfxy(175.00,60.00)}\pgflineto{\pgfxy(175.00,65.00)}\pgflineto{\pgfxy(170.00,65.00)}\pgfclosepath\pgffill
\color[rgb]{0,0,0}\pgfputat{\pgfxy(172.50,61.50)}{\pgfbox[bottom,left]{\fontsize{11.38}{13.66}\selectfont \makebox[0pt]{$v$}}}
\pgfmoveto{\pgfxy(160.00,50.00)}\pgflineto{\pgfxy(200.00,50.00)}\pgflineto{\pgfxy(200.00,55.00)}\pgflineto{\pgfxy(195.00,55.00)}\pgflineto{\pgfxy(195.00,65.00)}\pgflineto{\pgfxy(190.00,65.00)}\pgflineto{\pgfxy(190.00,70.00)}\pgflineto{\pgfxy(185.00,70.00)}\pgflineto{\pgfxy(185.00,75.00)}\pgflineto{\pgfxy(170.00,75.00)}\pgflineto{\pgfxy(170.00,85.00)}\pgflineto{\pgfxy(165.00,85.00)}\pgflineto{\pgfxy(165.00,90.00)}\pgflineto{\pgfxy(160.00,90.00)}\pgfclosepath\pgfstroke
\pgfsetlinewidth{0.60mm}\pgfmoveto{\pgfxy(170.00,75.00)}\pgflineto{\pgfxy(170.00,70.00)}\pgflineto{\pgfxy(180.00,70.00)}\pgflineto{\pgfxy(180.00,65.00)}\pgflineto{\pgfxy(185.00,65.00)}\pgflineto{\pgfxy(185.00,60.00)}\pgflineto{\pgfxy(195.00,60.00)}\pgflineto{\pgfxy(195.00,65.00)}\pgflineto{\pgfxy(190.00,65.00)}\pgflineto{\pgfxy(190.00,70.00)}\pgflineto{\pgfxy(185.00,70.00)}\pgflineto{\pgfxy(185.00,75.00)}\pgfclosepath\pgfstroke
\pgfputat{\pgfxy(185.00,66.50)}{\pgfbox[bottom,left]{\fontsize{11.38}{13.66}\selectfont \makebox[0pt]{$R_v$}}}
\pgfmoveto{\pgfxy(120.00,75.00)}\pgflineto{\pgfxy(120.00,60.00)}\pgflineto{\pgfxy(145.00,60.00)}\pgflineto{\pgfxy(145.00,65.00)}\pgflineto{\pgfxy(125.00,65.00)}\pgflineto{\pgfxy(125.00,75.00)}\pgfclosepath\pgfstroke
\pgfputat{\pgfxy(135.00,61.50)}{\pgfbox[bottom,left]{\fontsize{11.38}{13.66}\selectfont \makebox[0pt]{$H_v$}}}
\end{pgfpicture}%
$$
\caption{The hook and rim hook  of $(\lambda, v)$}
\label{fig:hookandrimhook}
\end{figure}

Let $\lambda$ be  a partition. Denote its {\em conjugate} by $\lambda'=(\lambda_1',\lambda_2', \ldots)$, where
$\lambda_i'$ is the number of parts of $\lambda$ that are $\geq i$.
For a cell $v\in \lambda$, the {\em hook} $H_v:=H_v(\lambda)$ of the cell $v$ in $\lambda$ is the set of all cells  of $\lambda$ lying in the same column above $v$ or in the same row to the right of $v$, including $v$ itself. Obviously, the number of all cells in a hook $H_v$ equals the hook length $h_v$ for any $v$ in $\lambda$.
A cell $v$ in $\lambda$ is called {\em boundary cell} if the upright corner of $v$ is in boundaries of $\lambda$.
A {\em border strip} is a sequence  $x_{0}, x_{1},\ldots, x_{m}$ of boundary cells in $\lambda$
such that $x_{j-1}$ and $x_{j}$ have a common side for $1\leq j\leq m$.
The {\em rim hook} $R_v:=R_v(\lambda)$ of $(\lambda,v)$  is a border strip $x_{0}, x_{1},\ldots, x_{m}$ of cells in $\lambda$ such that
$x_{0}$ (resp. $x_{m}$) is the uppermost  (resp.  rightmost) cell of the hook $H_v$.
The \emph{rim hook length} $r_v$ of $v$ in $\lambda$ is defined to be  the number of all cells in the  rim hook $R_v$.
It is easy to see that the hook $H_v$ and the rim hook $R_v$ have the same length, same height (or number of rows), and same width(or number of columns). See Figure~\ref{fig:hookandrimhook}.

For a given nonnegative integers $a$ and $m$,
let $\mathcal A$ be the set of partitions whose largest part is bounded by $m$,
$\tilde{\mathcal A}$ the set of partitions whose largest part is bounded by $m$ and parts are at most $a$,
and $\mathcal R$ the set of nondecreasing sequences $(r_1,\ldots, r_t)$, $t\ge 0$, with
\begin{equation*}\label{eq:increaing}
a+1 \leq r_1 \leq r_2 \leq \cdots \leq r_t \leq a+m.
\end{equation*}

\begin{figure}[t]
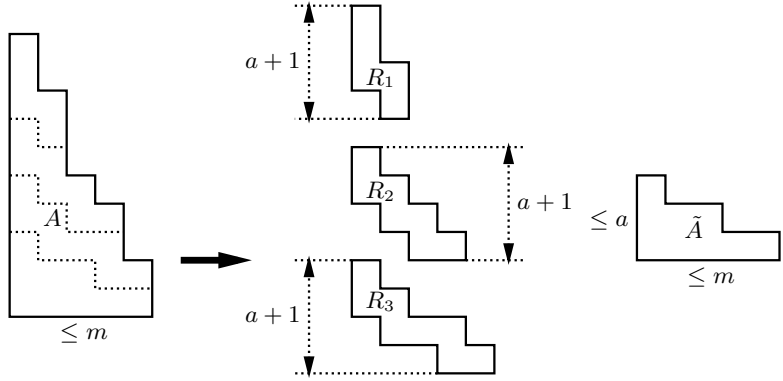

$$
\centering
\begin{pgfpicture}{50.50mm}{58.00mm}{155.75mm}{110.75mm}
\pgfsetxvec{\pgfpoint{0.75mm}{0mm}}
\pgfsetyvec{\pgfpoint{0mm}{0.75mm}}
\color[rgb]{0,0,0}\pgfsetlinewidth{0.30mm}\pgfsetdash{}{0mm}
\pgfputat{\pgfxy(77.50,106.00)}{\pgfbox[bottom,left]{\fontsize{8.54}{10.24}\selectfont \makebox[0pt]{$A$}}}
\pgfsetlinewidth{0.90mm}\pgfmoveto{\pgfxy(100.00,100.00)}\pgflineto{\pgfxy(110.00,100.00)}\pgfstroke
\pgfmoveto{\pgfxy(110.00,100.00)}\pgflineto{\pgfxy(107.20,100.70)}\pgflineto{\pgfxy(107.20,99.30)}\pgflineto{\pgfxy(110.00,100.00)}\pgfclosepath\pgffill
\pgfmoveto{\pgfxy(110.00,100.00)}\pgflineto{\pgfxy(107.20,100.70)}\pgflineto{\pgfxy(107.20,99.30)}\pgflineto{\pgfxy(110.00,100.00)}\pgfclosepath\pgfstroke
\pgfputat{\pgfxy(178.50,106.00)}{\pgfbox[bottom,left]{\fontsize{8.54}{10.24}\selectfont \makebox[0pt][r]{$\leq a$}}}
\pgfputat{\pgfxy(83.00,86.00)}{\pgfbox[bottom,left]{\fontsize{8.54}{10.24}\selectfont \makebox[0pt]{$\leq m$}}}
\pgfsetlinewidth{0.30mm}\pgfmoveto{\pgfxy(70.00,140.00)}\pgflineto{\pgfxy(70.00,90.00)}\pgflineto{\pgfxy(95.00,90.00)}\pgflineto{\pgfxy(95.00,100.00)}\pgflineto{\pgfxy(90.00,100.00)}\pgflineto{\pgfxy(90.00,110.00)}\pgflineto{\pgfxy(85.00,110.00)}\pgflineto{\pgfxy(85.00,115.00)}\pgflineto{\pgfxy(80.00,115.00)}\pgflineto{\pgfxy(80.00,130.00)}\pgflineto{\pgfxy(75.00,130.00)}\pgflineto{\pgfxy(75.00,140.00)}\pgfclosepath\pgfstroke
\pgfmoveto{\pgfxy(130.00,145.00)}\pgflineto{\pgfxy(130.00,130.00)}\pgflineto{\pgfxy(135.00,130.00)}\pgflineto{\pgfxy(135.00,125.00)}\pgflineto{\pgfxy(140.00,125.00)}\pgflineto{\pgfxy(140.00,135.00)}\pgflineto{\pgfxy(135.00,135.00)}\pgflineto{\pgfxy(135.00,145.00)}\pgfclosepath\pgfstroke
\pgfmoveto{\pgfxy(130.00,100.00)}\pgflineto{\pgfxy(130.00,90.00)}\pgflineto{\pgfxy(135.00,90.00)}\pgflineto{\pgfxy(135.00,85.00)}\pgflineto{\pgfxy(145.00,85.00)}\pgflineto{\pgfxy(145.00,80.00)}\pgflineto{\pgfxy(155.00,80.00)}\pgflineto{\pgfxy(155.00,85.00)}\pgflineto{\pgfxy(150.00,85.00)}\pgflineto{\pgfxy(150.00,90.00)}\pgflineto{\pgfxy(140.00,90.00)}\pgflineto{\pgfxy(140.00,95.00)}\pgflineto{\pgfxy(135.00,95.00)}\pgflineto{\pgfxy(135.00,100.00)}\pgfclosepath\pgfstroke
\pgfputat{\pgfxy(190.00,104.00)}{\pgfbox[bottom,left]{\fontsize{8.54}{10.24}\selectfont \makebox[0pt]{$\tilde{A}$}}}
\pgfsetdash{{0.30mm}{0.50mm}}{0mm}\pgfmoveto{\pgfxy(135.00,125.00)}\pgflineto{\pgfxy(120.00,125.00)}\pgfstroke
\pgfmoveto{\pgfxy(160.00,120.00)}\pgflineto{\pgfxy(135.00,120.00)}\pgfstroke
\pgfmoveto{\pgfxy(150.00,100.00)}\pgflineto{\pgfxy(160.00,100.00)}\pgfstroke
\pgfmoveto{\pgfxy(130.00,100.00)}\pgflineto{\pgfxy(120.00,100.00)}\pgfstroke
\pgfmoveto{\pgfxy(145.00,80.00)}\pgflineto{\pgfxy(120.00,80.00)}\pgfstroke
\pgfmoveto{\pgfxy(122.50,145.00)}\pgflineto{\pgfxy(122.50,125.00)}\pgfstroke
\pgfmoveto{\pgfxy(122.50,145.00)}\pgflineto{\pgfxy(121.80,142.20)}\pgflineto{\pgfxy(123.20,142.20)}\pgflineto{\pgfxy(122.50,145.00)}\pgfclosepath\pgffill
\pgfsetdash{}{0mm}\pgfmoveto{\pgfxy(122.50,145.00)}\pgflineto{\pgfxy(121.80,142.20)}\pgflineto{\pgfxy(123.20,142.20)}\pgflineto{\pgfxy(122.50,145.00)}\pgfclosepath\pgfstroke
\pgfmoveto{\pgfxy(122.50,125.00)}\pgflineto{\pgfxy(123.20,127.80)}\pgflineto{\pgfxy(121.80,127.80)}\pgflineto{\pgfxy(122.50,125.00)}\pgfclosepath\pgffill
\pgfmoveto{\pgfxy(122.50,125.00)}\pgflineto{\pgfxy(123.20,127.80)}\pgflineto{\pgfxy(121.80,127.80)}\pgflineto{\pgfxy(122.50,125.00)}\pgfclosepath\pgfstroke
\pgfsetdash{{0.30mm}{0.50mm}}{0mm}\pgfmoveto{\pgfxy(157.50,120.00)}\pgflineto{\pgfxy(157.50,100.00)}\pgfstroke
\pgfmoveto{\pgfxy(157.50,120.00)}\pgflineto{\pgfxy(156.80,117.20)}\pgflineto{\pgfxy(158.20,117.20)}\pgflineto{\pgfxy(157.50,120.00)}\pgfclosepath\pgffill
\pgfsetdash{}{0mm}\pgfmoveto{\pgfxy(157.50,120.00)}\pgflineto{\pgfxy(156.80,117.20)}\pgflineto{\pgfxy(158.20,117.20)}\pgflineto{\pgfxy(157.50,120.00)}\pgfclosepath\pgfstroke
\pgfmoveto{\pgfxy(157.50,100.00)}\pgflineto{\pgfxy(158.20,102.80)}\pgflineto{\pgfxy(156.80,102.80)}\pgflineto{\pgfxy(157.50,100.00)}\pgfclosepath\pgffill
\pgfmoveto{\pgfxy(157.50,100.00)}\pgflineto{\pgfxy(158.20,102.80)}\pgflineto{\pgfxy(156.80,102.80)}\pgflineto{\pgfxy(157.50,100.00)}\pgfclosepath\pgfstroke
\pgfsetdash{{0.30mm}{0.50mm}}{0mm}\pgfmoveto{\pgfxy(122.50,100.00)}\pgflineto{\pgfxy(122.50,80.00)}\pgfstroke
\pgfmoveto{\pgfxy(122.50,100.00)}\pgflineto{\pgfxy(121.80,97.20)}\pgflineto{\pgfxy(123.20,97.20)}\pgflineto{\pgfxy(122.50,100.00)}\pgfclosepath\pgffill
\pgfsetdash{}{0mm}\pgfmoveto{\pgfxy(122.50,100.00)}\pgflineto{\pgfxy(121.80,97.20)}\pgflineto{\pgfxy(123.20,97.20)}\pgflineto{\pgfxy(122.50,100.00)}\pgfclosepath\pgfstroke
\pgfmoveto{\pgfxy(122.50,80.00)}\pgflineto{\pgfxy(123.20,82.80)}\pgflineto{\pgfxy(121.80,82.80)}\pgflineto{\pgfxy(122.50,80.00)}\pgfclosepath\pgffill
\pgfmoveto{\pgfxy(122.50,80.00)}\pgflineto{\pgfxy(123.20,82.80)}\pgflineto{\pgfxy(121.80,82.80)}\pgflineto{\pgfxy(122.50,80.00)}\pgfclosepath\pgfstroke
\pgfputat{\pgfxy(121.00,134.00)}{\pgfbox[bottom,left]{\fontsize{8.54}{10.24}\selectfont \makebox[0pt][r]{$a+1$}}}
\pgfputat{\pgfxy(121.00,89.00)}{\pgfbox[bottom,left]{\fontsize{8.54}{10.24}\selectfont \makebox[0pt][r]{$a+1$}}}
\pgfputat{\pgfxy(159.00,109.00)}{\pgfbox[bottom,left]{\fontsize{8.54}{10.24}\selectfont $a+1$}}
\pgfputat{\pgfxy(135.00,131.00)}{\pgfbox[bottom,left]{\fontsize{8.54}{10.24}\selectfont \makebox[0pt]{$R_1$}}}
\pgfputat{\pgfxy(135.00,111.00)}{\pgfbox[bottom,left]{\fontsize{8.54}{10.24}\selectfont \makebox[0pt]{$R_2$}}}
\pgfputat{\pgfxy(135.00,91.50)}{\pgfbox[bottom,left]{\fontsize{8.54}{10.24}\selectfont \makebox[0pt]{$R_3$}}}
\pgfputat{\pgfxy(193.00,96.00)}{\pgfbox[bottom,left]{\fontsize{8.54}{10.24}\selectfont \makebox[0pt]{$\leq m$}}}
\pgfmoveto{\pgfxy(130.00,120.00)}\pgflineto{\pgfxy(130.00,110.00)}\pgflineto{\pgfxy(135.00,110.00)}\pgflineto{\pgfxy(135.00,105.00)}\pgflineto{\pgfxy(140.00,105.00)}\pgflineto{\pgfxy(140.00,100.00)}\pgflineto{\pgfxy(150.00,100.00)}\pgflineto{\pgfxy(150.00,105.00)}\pgflineto{\pgfxy(145.00,105.00)}\pgflineto{\pgfxy(145.00,110.00)}\pgflineto{\pgfxy(140.00,110.00)}\pgflineto{\pgfxy(140.00,115.00)}\pgflineto{\pgfxy(135.00,115.00)}\pgflineto{\pgfxy(135.00,120.00)}\pgfclosepath\pgfstroke
\pgfmoveto{\pgfxy(180.00,115.00)}\pgflineto{\pgfxy(180.00,100.00)}\pgflineto{\pgfxy(205.00,100.00)}\pgflineto{\pgfxy(205.00,105.00)}\pgflineto{\pgfxy(195.00,105.00)}\pgflineto{\pgfxy(195.00,110.00)}\pgflineto{\pgfxy(185.00,110.00)}\pgflineto{\pgfxy(185.00,115.00)}\pgfclosepath\pgfstroke
\pgfsetdash{{0.30mm}{0.50mm}}{0mm}\pgfmoveto{\pgfxy(70.00,125.00)}\pgflineto{\pgfxy(75.00,125.00)}\pgflineto{\pgfxy(75.00,120.00)}\pgflineto{\pgfxy(80.00,120.00)}\pgfstroke
\pgfmoveto{\pgfxy(70.00,115.00)}\pgflineto{\pgfxy(75.00,115.00)}\pgflineto{\pgfxy(75.00,110.00)}\pgflineto{\pgfxy(80.00,110.00)}\pgflineto{\pgfxy(80.00,105.00)}\pgflineto{\pgfxy(90.00,105.00)}\pgfstroke
\pgfmoveto{\pgfxy(70.00,105.00)}\pgflineto{\pgfxy(75.00,105.00)}\pgflineto{\pgfxy(75.00,100.00)}\pgflineto{\pgfxy(85.00,100.00)}\pgflineto{\pgfxy(85.00,95.00)}\pgflineto{\pgfxy(95.00,95.00)}\pgfstroke
\pgfmoveto{\pgfxy(130.00,145.00)}\pgflineto{\pgfxy(120.00,145.00)}\pgfstroke
\end{pgfpicture}%
$$
\caption{Move out rim hooks from $A$ to get $(\tilde{A}; r_1, \ldots, r_t)$}
\label{fig:rimhook}
\end{figure}

We first describe an important algorithm, called {\em Pealing Algorithm}, which will be used in the construction of the involution $\Phi$.
From the leftmost top cell of the diagram of $A \in \mathcal A$,
remove a rim hook of height $a+1$, if any, in such a way that what remains is a diagram of a partition, and continue removing rim hooks of height $a+1$ in this way as long as possible.
Denote the remained partition by $\tilde{A}$ as shown in Figure~\ref{fig:rimhook}.
Clearly, the length of $\tilde{A}$ is less than or equal to $a$ and $\tilde{A} \in \tilde{\mathcal A}$.

Each removed rim hook in the above  transformation corresponds to some cell in the first column of $A$.
Let  $v_1, \ldots, v_t$ be the  cells from top to bottom
 corresponding to removed rim hooks.
If $r_1, \ldots, r_t$ are the lengths of removed rim hooks $R_1, \ldots, R_t$, then
$$
r_i=a_{v_i}+a+1 \quad \text{for all $1\le i\le t$}.
$$
Since $A$  is a partition with largest part bounded by $m$, we have $$0\le a_{v_{1}}\leq a_{v_{2}}\leq \cdots \leq a_{v_{t}}\le m-1$$ and $(r_1, \ldots, r_t)\in \mathcal R$.

\begin{lem}[Pealing Algorithm]\label{lem:pealing}
The mapping $A \mapsto (\tilde{A}; r_1, \ldots, r_t)$ is a bijection from $\mathcal A$ to $\tilde{\mathcal A} \times \mathcal R$.
\end{lem}

\begin{proof}
It is sufficient to construct the inverse of pealing algorithm.

Starting from $(\tilde{A}; r_1, \ldots, r_t) \in \tilde{\mathcal A} \times \mathcal R$, let $A_t:=\tilde{A}=(\lambda_1,\ldots,\lambda_{\ell})$ with $\ell \le a$.
If there exists a cell whose arm length $< r_t-a-1$, let $v_t:=(1,\alpha)$ be the lowest such cell. Otherwise, let $v_t:=(1,\alpha)$ with $\alpha=\ell(A_t)+1$. Since the length of $A_t$ is less than or equal to $a$, we have $l_{v_t}(A_t) < a$.
Define the partition
$$A_{t-1}:=(\lambda_1, \ldots, \lambda_{\alpha-1}, r_t-a, \lambda_{\alpha}+1, \ldots, \lambda_{\alpha+a-1}+1),$$
where $\lambda_i :=0$ for $i>\ell$.
Clearly
\begin{equation}\label{eq:rimhook1}
l_{v_t}(A_{t-1})=a \quad\text{and}\quad a_{v_t}(A_{t-1})=r_t-a-1.
\end{equation}

%

Next, we proceed by induction on $i$ from $t-1$ to $1$.
Suppose that we have found $A_{t-1}, A_{t-2}, \ldots, A_{i} = (\lambda_1,\ldots,\lambda_\ell)$ with
\begin{equation}\label{eq:rimhook2}
l_{v_{i+1}}(A_{i})=a \quad\text{and}\quad a_{v_{i+1}}(A_{i})=r_{i+1}-a-1.
\end{equation}
If there exists a cell whose arm length $< r_i-a-1$, let $v_i:=(1,\alpha)$ be the lowest such cell. Otherwise, let $v_i:=(1,\alpha)$ with $\alpha=\ell(A_i)+1$.
Since we have $a_{v_{i+1}}(A_{i})=r_{i+1}-a-1\ge r_{i}-a-1$ by \eqref{eq:rimhook2}, the cell $v_i$ is above $v_{i+1}$. Also, since $l_{v_{i+1}}(A_{i}) = a$ by \eqref{eq:rimhook2}, we have $l_{v_{i}}(A_{i}) < a$.
Let $$A_{i-1}:=(\lambda_1, \ldots, \lambda_{\alpha-1}, r_i-a, \lambda_{\alpha}+1, \ldots, \lambda_{\alpha+a-1}+1),$$
where $\lambda_i :=0$ for $i>\ell$. Clearly,
\begin{equation}\label{eq:rimhook3}
l_{v_{i}}(A_{i-1})=a \quad\text{and}\quad a_{v_{i}}(A_{i-1})=r_{i}-a-1.
\end{equation}


Finally, define $A:=A_0 \in \mathcal{A}$. Applying pealing algorithm to such a partition $A$ induced from $(\tilde{A}; r_1, \ldots, r_t)$, by \eqref{eq:rimhook1} and \eqref{eq:rimhook3}, we get $(\tilde{A}; r_1, \ldots, r_t)$ again.

\end{proof}

For nonnegative integers $m$ and $n$, the $q$-ascending factorial is defined by
$(a;q)_n = (1-a) (1-aq) \cdots (1-aq^{n-1})$ and the $q$-binomial coefficient is defined by
\begin{eqnarray*}
\qbin{n}{m} = \frac{(q;q)_n}{(q;q)_m (q;q)_{n-m}} \quad \text{for $0\leq m \leq n$}.
\end{eqnarray*}

\rmk The pealing algorithm gives a bijective proof of the formula
\begin{align*}
\frac{1}{(q;q)_m} &= \qbin{m+a}{a} \times \frac{1}{(q^{a+1};q)_m}\\
&=\qbin{m+a}{a} \times \sum_{t\ge 0} q^{t(a+1)} \qbin{m-1+t}{t},
\end{align*}
where $t$ means the number of removed rim hooks. See \cite[Eq. (3.3.7)]{MR1634067}.

\section{Main result}
Let $\mathcal{F}_n(a,l,m)$ be the set of pointed partitions $(\lambda,v)$ of $n$ such that  $a_v=a$, $l_v=l$ and $m_v=m$. We shall divide the construction of the involution $\Phi$ in three basic steps.

\subsection*{Step 1: Bessenrodt-Han decomposition $\varphi_{a,l,m}: (\lambda,v)\mapsto (A,B,C,D,E)$}
Let $\mathcal{Q}_n(a,l,m)$ be the set of quintuples $(A,B,C,D,E)$ such that
$A$ is a partition whose largest part is bounded by  $m$,
$B$ is a partition whose diagram fits inside an $l\times a$ rectangle,
$C$ is a partition whose all parts are greater than or equal to $m+a+1$,
$D$ is a partition whose diagram is an $(l+1) \times (m+1)$ rectangle,
$E$ is a partition whose diagram is an $1 \times a$ rectangle,
and $$\abs{A}+\abs{B}+\abs{C}+\abs{D}+\abs{E}=n.$$

\begin{figure}[t]
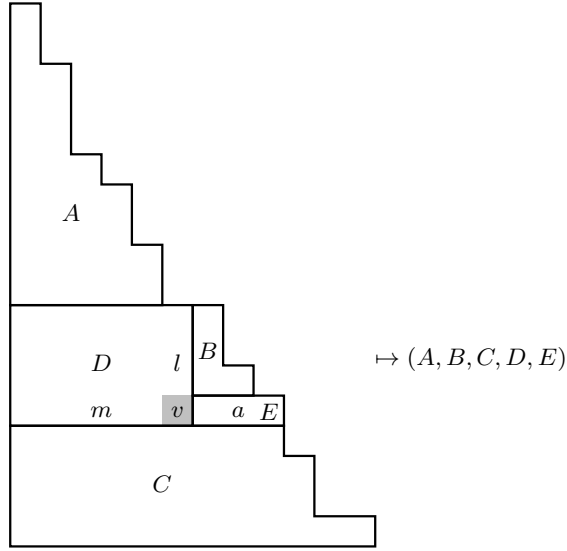

$$
\centering
\begin{pgfpicture}{-22.00mm}{30.00mm}{57.11mm}{106.00mm}
\pgfsetxvec{\pgfpoint{0.80mm}{0mm}}
\pgfsetyvec{\pgfpoint{0mm}{0.80mm}}
\color[rgb]{0,0,0}\pgfsetlinewidth{0.30mm}\pgfsetdash{}{0mm}
\color[rgb]{0.75294,0.75294,0.75294}\pgfmoveto{\pgfxy(0.00,60.00)}\pgflineto{\pgfxy(5.00,60.00)}\pgflineto{\pgfxy(5.00,65.00)}\pgflineto{\pgfxy(0.00,65.00)}\pgfclosepath\pgffill
\color[rgb]{0,0,0}\pgfmoveto{\pgfxy(-25.00,60.00)}\pgflineto{\pgfxy(-25.00,40.00)}\pgflineto{\pgfxy(35.00,40.00)}\pgflineto{\pgfxy(35.00,45.00)}\pgflineto{\pgfxy(25.00,45.00)}\pgflineto{\pgfxy(25.00,55.00)}\pgflineto{\pgfxy(20.00,55.00)}\pgflineto{\pgfxy(20.00,60.00)}\pgfclosepath\pgfstroke
\pgfmoveto{\pgfxy(5.00,80.00)}\pgflineto{\pgfxy(5.00,65.00)}\pgflineto{\pgfxy(15.00,65.00)}\pgflineto{\pgfxy(15.00,70.00)}\pgflineto{\pgfxy(10.00,70.00)}\pgflineto{\pgfxy(10.00,80.00)}\pgfclosepath\pgfstroke
\pgfputat{\pgfxy(-10.00,69.00)}{\pgfbox[bottom,left]{\fontsize{9.10}{10.93}\selectfont \makebox[0pt]{$D$}}}
\pgfputat{\pgfxy(0.00,49.00)}{\pgfbox[bottom,left]{\fontsize{9.10}{10.93}\selectfont \makebox[0pt]{$C$}}}
\pgfputat{\pgfxy(-15.00,94.00)}{\pgfbox[bottom,left]{\fontsize{9.10}{10.93}\selectfont \makebox[0pt]{$A$}}}
\pgfputat{\pgfxy(7.50,71.00)}{\pgfbox[bottom,left]{\fontsize{9.10}{10.93}\selectfont \makebox[0pt]{$B$}}}
\pgfputat{\pgfxy(2.50,61.50)}{\pgfbox[bottom,left]{\fontsize{9.10}{10.93}\selectfont \makebox[0pt]{$v$}}}
\pgfputat{\pgfxy(-10.00,61.50)}{\pgfbox[bottom,left]{\fontsize{9.10}{10.93}\selectfont \makebox[0pt]{$m$}}}
\pgfputat{\pgfxy(12.50,61.50)}{\pgfbox[bottom,left]{\fontsize{9.10}{10.93}\selectfont \makebox[0pt]{$a$}}}
\pgfputat{\pgfxy(2.50,69.00)}{\pgfbox[bottom,left]{\fontsize{9.10}{10.93}\selectfont \makebox[0pt]{$l$}}}
\pgfmoveto{\pgfxy(-25.00,60.00)}\pgflineto{\pgfxy(5.00,60.00)}\pgflineto{\pgfxy(5.00,80.00)}\pgflineto{\pgfxy(-25.00,80.00)}\pgfclosepath\pgfstroke
\pgfmoveto{\pgfxy(5.00,60.00)}\pgflineto{\pgfxy(20.00,60.00)}\pgflineto{\pgfxy(20.00,65.00)}\pgflineto{\pgfxy(5.00,65.00)}\pgfclosepath\pgfstroke
\pgfputat{\pgfxy(17.50,61.00)}{\pgfbox[bottom,left]{\fontsize{9.10}{10.93}\selectfont \makebox[0pt]{$E$}}}
\pgfmoveto{\pgfxy(-25.00,130.00)}\pgflineto{\pgfxy(-25.00,80.00)}\pgflineto{\pgfxy(0.00,80.00)}\pgflineto{\pgfxy(0.00,90.00)}\pgflineto{\pgfxy(-5.00,90.00)}\pgflineto{\pgfxy(-5.00,100.00)}\pgflineto{\pgfxy(-10.00,100.00)}\pgflineto{\pgfxy(-10.00,105.00)}\pgflineto{\pgfxy(-15.00,105.00)}\pgflineto{\pgfxy(-15.00,120.00)}\pgflineto{\pgfxy(-20.00,120.00)}\pgflineto{\pgfxy(-20.00,130.00)}\pgfclosepath\pgfstroke
\pgfputat{\pgfxy(35.00,70.00)}{\pgfbox[bottom,left]{\fontsize{9.10}{10.93}\selectfont $\mapsto (A,B,C,D,E)$}}
\end{pgfpicture}%
$$
\caption{Bessenrodt-Han decomposition of a pointed partition $(\lambda, v)$ }
\label{fig:decomposition}
\end{figure}

Suppose that $(\lambda,v)\in \mathcal{F}_n(a,l,m)$.
First of all, we decompose the diagram of $\lambda$ into five regions $(A,B,C,D,E)$ as shown in Figure~\ref{fig:decomposition}, modified slightly from the decomposition made in \cite{BESSENRODT:2009:HAL-00395700:1}, in which $E$ is implicit.
Denote this decomposition by $\varphi_{a,l,m}$.

Merging the other partitions $B$, $C$, $D$, and $E$ in bottom of $A$, we recover a pointed partition $(\lambda, v)$. So $\varphi_{a,l,m}$ is a bijection from $\mathcal{F}_n(a,l,m)$ to $\mathcal{Q}_n(a,l,m)$.

\begin{ex}
For the pointed partition $(\lambda, v)\in \mathcal{F}_{101}(3,3,5)$ with
  $$\lambda=(12,10,10,9,9,8,7,7,5,5,4,4,3,2,2,2,1,1)\quad \textrm{and}\quad  v=(6,5),
  $$
the decomposition
is given by
$(A,B,C,D,E)$ with $A=(5,5,4,4,3,2,2,2,1,1)$, $B=(2,1,1)$, $C=(12,10,10,9)$, $D=(6,6,6,6)$ and $E=(3)$.
\end{ex}

\subsection*{Step 2:  Transformation $\psi_{a,l,m}:  (A, B,C,D,E)\mapsto (\tilde{A}, B,\tilde{C},D,E)$}\label{step}

Let $\widetilde\mathcal{Q}_n(a,l,m)$ be the set of quintuples $(\tilde{A},B,\tilde{C},D,E)$ such that
$\tilde{A}$ is a partition whose diagram fits inside an $a\times m$ rectangle,
$B$ is a partition whose diagram fits inside an $l\times a$ rectangle,
$\tilde{C}$ is a partition whose all parts are greater than or equal to $a+1$,
$D$ is a partition whose diagram is an $(l+1) \times (m+1)$ rectangle,
$E$ is a partition whose diagram is an $1 \times a$ rectangle,
and $$\tilde{\abs{A}}+\abs{B}+\tilde{\abs{C}}+\abs{D}+\abs{E}=n.$$
Note that the partitions $\tilde{A}$, $B$, $\tilde{C}$, and $E$ could be empty.
The mapping $\psi_{a,l,m}$ from $\mathcal{Q}_n(a,l,m)$ to $\widetilde\mathcal{Q}_n(a,l,m)$ is defined by $\psi_{a,l,m}(A,B,C,D,E)=(\tilde{A},B,\tilde{C},D,E)$ as follows:
\begin{itemize}
\item Applying the pealing algorithm to $A$, we have $(\tilde{A}; r_1, \ldots, r_t)$.
In the above example, we have $\tilde{A}=(5,3,1)$ and $(r_1,r_2,r_3)=(5,7,8)$ with $t=3$.

\item Gluing $r_t, \ldots, r_1$ to $C$, we get the partition $\tilde{C}=(C, r_t, r_{t-1}, \ldots, r_1)$ obtained by adding $r_i$'s to the partition $C$. Clearly, all parts of $\tilde{C}$ are greater than or equal to $a+1$.
\end{itemize}


\begin{figure}[t]
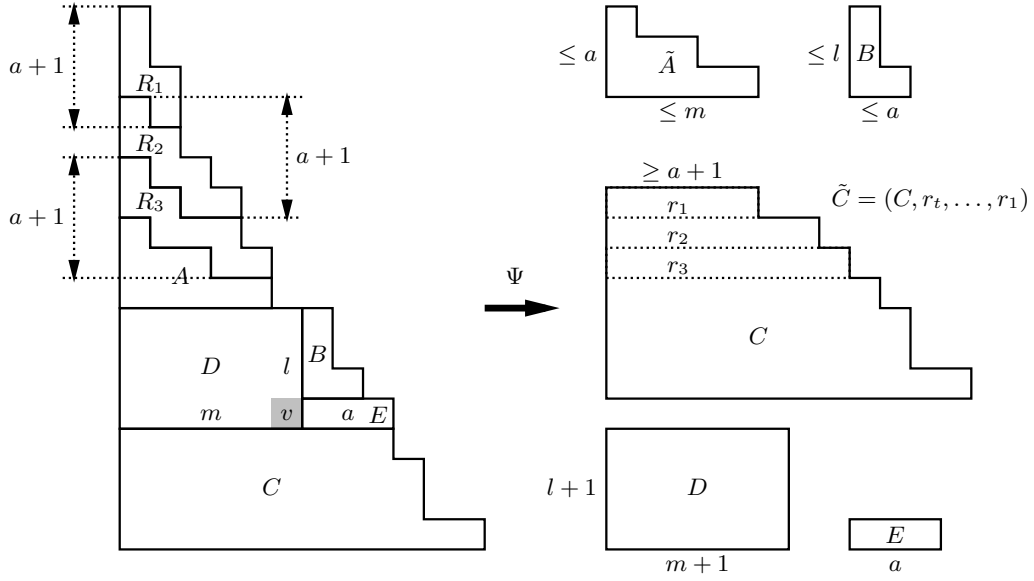

$$
\centering
\begin{pgfpicture}{25.97mm}{26.51mm}{154.00mm}{106.00mm}
\pgfsetxvec{\pgfpoint{0.80mm}{0mm}}
\pgfsetyvec{\pgfpoint{0mm}{0.80mm}}
\color[rgb]{0,0,0}\pgfsetlinewidth{0.30mm}\pgfsetdash{}{0mm}
\color[rgb]{0.75294,0.75294,0.75294}\pgfmoveto{\pgfxy(75.00,60.00)}\pgflineto{\pgfxy(80.00,60.00)}\pgflineto{\pgfxy(80.00,65.00)}\pgflineto{\pgfxy(75.00,65.00)}\pgfclosepath\pgffill
\color[rgb]{0,0,0}\pgfmoveto{\pgfxy(50.00,60.00)}\pgflineto{\pgfxy(50.00,40.00)}\pgflineto{\pgfxy(110.00,40.00)}\pgflineto{\pgfxy(110.00,45.00)}\pgflineto{\pgfxy(100.00,45.00)}\pgflineto{\pgfxy(100.00,55.00)}\pgflineto{\pgfxy(95.00,55.00)}\pgflineto{\pgfxy(95.00,60.00)}\pgfclosepath\pgfstroke
\pgfmoveto{\pgfxy(50.00,130.00)}\pgflineto{\pgfxy(50.00,115.00)}\pgflineto{\pgfxy(55.00,115.00)}\pgflineto{\pgfxy(55.00,110.00)}\pgflineto{\pgfxy(60.00,110.00)}\pgflineto{\pgfxy(60.00,120.00)}\pgflineto{\pgfxy(55.00,120.00)}\pgflineto{\pgfxy(55.00,130.00)}\pgfclosepath\pgfstroke
\pgfmoveto{\pgfxy(50.00,115.00)}\pgflineto{\pgfxy(50.00,105.00)}\pgflineto{\pgfxy(55.00,105.00)}\pgflineto{\pgfxy(55.00,100.00)}\pgflineto{\pgfxy(60.00,100.00)}\pgflineto{\pgfxy(60.00,95.00)}\pgflineto{\pgfxy(65.00,95.00)}\pgflineto{\pgfxy(70.00,95.00)}\pgflineto{\pgfxy(70.00,100.00)}\pgflineto{\pgfxy(65.00,100.00)}\pgflineto{\pgfxy(65.00,105.00)}\pgflineto{\pgfxy(60.00,105.00)}\pgflineto{\pgfxy(60.00,110.00)}\pgflineto{\pgfxy(55.00,110.00)}\pgflineto{\pgfxy(55.00,115.00)}\pgfclosepath\pgfstroke
\pgfmoveto{\pgfxy(50.00,105.00)}\pgflineto{\pgfxy(50.00,95.00)}\pgflineto{\pgfxy(55.00,95.00)}\pgflineto{\pgfxy(55.00,90.00)}\pgflineto{\pgfxy(65.00,90.00)}\pgflineto{\pgfxy(65.00,85.00)}\pgflineto{\pgfxy(75.00,85.00)}\pgflineto{\pgfxy(75.00,90.00)}\pgflineto{\pgfxy(70.00,90.00)}\pgflineto{\pgfxy(70.00,95.00)}\pgflineto{\pgfxy(60.00,95.00)}\pgflineto{\pgfxy(60.00,100.00)}\pgflineto{\pgfxy(55.00,100.00)}\pgflineto{\pgfxy(55.00,105.00)}\pgfclosepath\pgfstroke
\pgfmoveto{\pgfxy(50.00,95.00)}\pgflineto{\pgfxy(50.00,80.00)}\pgflineto{\pgfxy(75.00,80.00)}\pgflineto{\pgfxy(75.00,85.00)}\pgflineto{\pgfxy(65.00,85.00)}\pgflineto{\pgfxy(65.00,90.00)}\pgflineto{\pgfxy(55.00,90.00)}\pgflineto{\pgfxy(55.00,95.00)}\pgfclosepath\pgfstroke
\pgfmoveto{\pgfxy(80.00,80.00)}\pgflineto{\pgfxy(80.00,65.00)}\pgflineto{\pgfxy(90.00,65.00)}\pgflineto{\pgfxy(90.00,70.00)}\pgflineto{\pgfxy(85.00,70.00)}\pgflineto{\pgfxy(85.00,80.00)}\pgfclosepath\pgfstroke
\pgfputat{\pgfxy(65.00,69.00)}{\pgfbox[bottom,left]{\fontsize{9.10}{10.93}\selectfont \makebox[0pt]{$D$}}}
\pgfputat{\pgfxy(75.00,49.00)}{\pgfbox[bottom,left]{\fontsize{9.10}{10.93}\selectfont \makebox[0pt]{$C$}}}
\pgfputat{\pgfxy(60.00,84.00)}{\pgfbox[bottom,left]{\fontsize{9.10}{10.93}\selectfont \makebox[0pt]{$A$}}}
\pgfputat{\pgfxy(82.50,71.00)}{\pgfbox[bottom,left]{\fontsize{9.10}{10.93}\selectfont \makebox[0pt]{$B$}}}
\pgfputat{\pgfxy(77.50,61.50)}{\pgfbox[bottom,left]{\fontsize{9.10}{10.93}\selectfont \makebox[0pt]{$v$}}}
\pgfputat{\pgfxy(65.00,61.50)}{\pgfbox[bottom,left]{\fontsize{9.10}{10.93}\selectfont \makebox[0pt]{$m$}}}
\pgfputat{\pgfxy(87.50,61.50)}{\pgfbox[bottom,left]{\fontsize{9.10}{10.93}\selectfont \makebox[0pt]{$a$}}}
\pgfputat{\pgfxy(77.50,69.00)}{\pgfbox[bottom,left]{\fontsize{9.10}{10.93}\selectfont \makebox[0pt]{$l$}}}
\pgfsetdash{{0.30mm}{0.50mm}}{0mm}\pgfmoveto{\pgfxy(50.00,130.00)}\pgflineto{\pgfxy(40.00,130.00)}\pgfstroke
\pgfmoveto{\pgfxy(55.00,110.00)}\pgflineto{\pgfxy(40.00,110.00)}\pgfstroke
\pgfmoveto{\pgfxy(80.00,115.00)}\pgflineto{\pgfxy(55.00,115.00)}\pgfstroke
\pgfmoveto{\pgfxy(70.00,95.00)}\pgflineto{\pgfxy(80.00,95.00)}\pgfstroke
\pgfmoveto{\pgfxy(50.00,105.00)}\pgflineto{\pgfxy(40.00,105.00)}\pgfstroke
\pgfmoveto{\pgfxy(65.00,85.00)}\pgflineto{\pgfxy(40.00,85.00)}\pgfstroke
\pgfmoveto{\pgfxy(42.50,130.00)}\pgflineto{\pgfxy(42.50,110.00)}\pgfstroke
\pgfmoveto{\pgfxy(42.50,130.00)}\pgflineto{\pgfxy(41.80,127.20)}\pgflineto{\pgfxy(43.20,127.20)}\pgflineto{\pgfxy(42.50,130.00)}\pgfclosepath\pgffill
\pgfsetdash{}{0mm}\pgfmoveto{\pgfxy(42.50,130.00)}\pgflineto{\pgfxy(41.80,127.20)}\pgflineto{\pgfxy(43.20,127.20)}\pgflineto{\pgfxy(42.50,130.00)}\pgfclosepath\pgfstroke
\pgfmoveto{\pgfxy(42.50,110.00)}\pgflineto{\pgfxy(43.20,112.80)}\pgflineto{\pgfxy(41.80,112.80)}\pgflineto{\pgfxy(42.50,110.00)}\pgfclosepath\pgffill
\pgfmoveto{\pgfxy(42.50,110.00)}\pgflineto{\pgfxy(43.20,112.80)}\pgflineto{\pgfxy(41.80,112.80)}\pgflineto{\pgfxy(42.50,110.00)}\pgfclosepath\pgfstroke
\pgfsetdash{{0.30mm}{0.50mm}}{0mm}\pgfmoveto{\pgfxy(77.50,115.00)}\pgflineto{\pgfxy(77.50,95.00)}\pgfstroke
\pgfmoveto{\pgfxy(77.50,115.00)}\pgflineto{\pgfxy(76.80,112.20)}\pgflineto{\pgfxy(78.20,112.20)}\pgflineto{\pgfxy(77.50,115.00)}\pgfclosepath\pgffill
\pgfsetdash{}{0mm}\pgfmoveto{\pgfxy(77.50,115.00)}\pgflineto{\pgfxy(76.80,112.20)}\pgflineto{\pgfxy(78.20,112.20)}\pgflineto{\pgfxy(77.50,115.00)}\pgfclosepath\pgfstroke
\pgfmoveto{\pgfxy(77.50,95.00)}\pgflineto{\pgfxy(78.20,97.80)}\pgflineto{\pgfxy(76.80,97.80)}\pgflineto{\pgfxy(77.50,95.00)}\pgfclosepath\pgffill
\pgfmoveto{\pgfxy(77.50,95.00)}\pgflineto{\pgfxy(78.20,97.80)}\pgflineto{\pgfxy(76.80,97.80)}\pgflineto{\pgfxy(77.50,95.00)}\pgfclosepath\pgfstroke
\pgfsetdash{{0.30mm}{0.50mm}}{0mm}\pgfmoveto{\pgfxy(42.50,105.00)}\pgflineto{\pgfxy(42.50,85.00)}\pgfstroke
\pgfmoveto{\pgfxy(42.50,105.00)}\pgflineto{\pgfxy(41.80,102.20)}\pgflineto{\pgfxy(43.20,102.20)}\pgflineto{\pgfxy(42.50,105.00)}\pgfclosepath\pgffill
\pgfsetdash{}{0mm}\pgfmoveto{\pgfxy(42.50,105.00)}\pgflineto{\pgfxy(41.80,102.20)}\pgflineto{\pgfxy(43.20,102.20)}\pgflineto{\pgfxy(42.50,105.00)}\pgfclosepath\pgfstroke
\pgfmoveto{\pgfxy(42.50,85.00)}\pgflineto{\pgfxy(43.20,87.80)}\pgflineto{\pgfxy(41.80,87.80)}\pgflineto{\pgfxy(42.50,85.00)}\pgfclosepath\pgffill
\pgfmoveto{\pgfxy(42.50,85.00)}\pgflineto{\pgfxy(43.20,87.80)}\pgflineto{\pgfxy(41.80,87.80)}\pgflineto{\pgfxy(42.50,85.00)}\pgfclosepath\pgfstroke
\pgfputat{\pgfxy(41.00,119.00)}{\pgfbox[bottom,left]{\fontsize{9.10}{10.93}\selectfont \makebox[0pt][r]{$a+1$}}}
\pgfputat{\pgfxy(41.00,94.00)}{\pgfbox[bottom,left]{\fontsize{9.10}{10.93}\selectfont \makebox[0pt][r]{$a+1$}}}
\pgfputat{\pgfxy(79.00,104.00)}{\pgfbox[bottom,left]{\fontsize{9.10}{10.93}\selectfont $a+1$}}
\pgfmoveto{\pgfxy(155.00,115.00)}\pgflineto{\pgfxy(130.00,115.00)}\pgflineto{\pgfxy(130.00,130.00)}\pgflineto{\pgfxy(135.00,130.00)}\pgflineto{\pgfxy(135.00,125.00)}\pgflineto{\pgfxy(145.00,125.00)}\pgflineto{\pgfxy(145.00,120.00)}\pgflineto{\pgfxy(155.00,120.00)}\pgfclosepath\pgfstroke
\pgfmoveto{\pgfxy(170.00,130.00)}\pgflineto{\pgfxy(170.00,115.00)}\pgflineto{\pgfxy(180.00,115.00)}\pgflineto{\pgfxy(180.00,120.00)}\pgflineto{\pgfxy(175.00,120.00)}\pgflineto{\pgfxy(175.00,130.00)}\pgfclosepath\pgfstroke
\pgfputat{\pgfxy(145.00,49.00)}{\pgfbox[bottom,left]{\fontsize{9.10}{10.93}\selectfont \makebox[0pt]{$D$}}}
\pgfputat{\pgfxy(140.00,119.00)}{\pgfbox[bottom,left]{\fontsize{9.10}{10.93}\selectfont \makebox[0pt]{$\tilde{A}$}}}
\pgfputat{\pgfxy(171.00,121.00)}{\pgfbox[bottom,left]{\fontsize{9.10}{10.93}\selectfont $B$}}
\pgfputat{\pgfxy(145.00,36.50)}{\pgfbox[bottom,left]{\fontsize{9.10}{10.93}\selectfont \makebox[0pt]{$m+1$}}}
\pgfputat{\pgfxy(177.50,36.50)}{\pgfbox[bottom,left]{\fontsize{9.10}{10.93}\selectfont \makebox[0pt]{$a$}}}
\pgfputat{\pgfxy(128.50,49.00)}{\pgfbox[bottom,left]{\fontsize{9.10}{10.93}\selectfont \makebox[0pt][r]{$l+1$}}}
\pgfputat{\pgfxy(55.00,116.00)}{\pgfbox[bottom,left]{\fontsize{9.10}{10.93}\selectfont \makebox[0pt]{$R_1$}}}
\pgfputat{\pgfxy(55.00,106.00)}{\pgfbox[bottom,left]{\fontsize{9.10}{10.93}\selectfont \makebox[0pt]{$R_2$}}}
\pgfputat{\pgfxy(55.00,96.50)}{\pgfbox[bottom,left]{\fontsize{9.10}{10.93}\selectfont \makebox[0pt]{$R_3$}}}
\pgfsetdash{{0.30mm}{0.50mm}}{0mm}\pgfmoveto{\pgfxy(130.00,85.00)}\pgflineto{\pgfxy(170.00,85.00)}\pgflineto{\pgfxy(170.00,90.00)}\pgflineto{\pgfxy(130.00,90.00)}\pgfclosepath\pgfstroke
\pgfmoveto{\pgfxy(130.00,95.00)}\pgflineto{\pgfxy(155.00,95.00)}\pgflineto{\pgfxy(155.00,100.00)}\pgflineto{\pgfxy(130.00,100.00)}\pgfclosepath\pgfstroke
\pgfputat{\pgfxy(140.00,96.00)}{\pgfbox[bottom,left]{\fontsize{9.10}{10.93}\selectfont $r_1$}}
\pgfputat{\pgfxy(140.00,91.00)}{\pgfbox[bottom,left]{\fontsize{9.10}{10.93}\selectfont $r_2$}}
\pgfputat{\pgfxy(140.00,86.00)}{\pgfbox[bottom,left]{\fontsize{9.10}{10.93}\selectfont $r_3$}}
\pgfsetdash{}{0mm}\pgfsetlinewidth{0.90mm}\pgfmoveto{\pgfxy(110.00,80.00)}\pgflineto{\pgfxy(120.00,80.00)}\pgfstroke
\pgfmoveto{\pgfxy(120.00,80.00)}\pgflineto{\pgfxy(117.20,80.70)}\pgflineto{\pgfxy(117.20,79.30)}\pgflineto{\pgfxy(120.00,80.00)}\pgfclosepath\pgffill
\pgfmoveto{\pgfxy(120.00,80.00)}\pgflineto{\pgfxy(117.20,80.70)}\pgflineto{\pgfxy(117.20,79.30)}\pgflineto{\pgfxy(120.00,80.00)}\pgfclosepath\pgfstroke
\pgfsetlinewidth{0.30mm}\pgfmoveto{\pgfxy(130.00,100.00)}\pgflineto{\pgfxy(130.00,65.00)}\pgflineto{\pgfxy(190.00,65.00)}\pgflineto{\pgfxy(190.00,70.00)}\pgflineto{\pgfxy(180.00,70.00)}\pgflineto{\pgfxy(180.00,80.00)}\pgflineto{\pgfxy(175.00,80.00)}\pgflineto{\pgfxy(175.00,85.00)}\pgflineto{\pgfxy(170.00,85.00)}\pgflineto{\pgfxy(170.00,90.00)}\pgflineto{\pgfxy(165.00,90.00)}\pgflineto{\pgfxy(165.00,95.00)}\pgflineto{\pgfxy(155.00,95.00)}\pgflineto{\pgfxy(155.00,100.00)}\pgfclosepath\pgfstroke
\pgfputat{\pgfxy(155.00,74.00)}{\pgfbox[bottom,left]{\fontsize{9.10}{10.93}\selectfont \makebox[0pt]{$C$}}}
\pgfputat{\pgfxy(115.00,84.00)}{\pgfbox[bottom,left]{\fontsize{9.10}{10.93}\selectfont \makebox[0pt]{$\Psi$}}}
\pgfputat{\pgfxy(175.00,111.50)}{\pgfbox[bottom,left]{\fontsize{9.10}{10.93}\selectfont \makebox[0pt]{$\leq a$}}}
\pgfputat{\pgfxy(128.50,121.00)}{\pgfbox[bottom,left]{\fontsize{9.10}{10.93}\selectfont \makebox[0pt][r]{$\leq a$}}}
\pgfputat{\pgfxy(142.50,111.50)}{\pgfbox[bottom,left]{\fontsize{9.10}{10.93}\selectfont \makebox[0pt]{$\leq m$}}}
\pgfputat{\pgfxy(168.50,121.00)}{\pgfbox[bottom,left]{\fontsize{9.10}{10.93}\selectfont \makebox[0pt][r]{$\leq l$}}}
\pgfputat{\pgfxy(142.50,101.50)}{\pgfbox[bottom,left]{\fontsize{9.10}{10.93}\selectfont \makebox[0pt]{$\geq a+1$}}}
\pgfmoveto{\pgfxy(50.00,60.00)}\pgflineto{\pgfxy(80.00,60.00)}\pgflineto{\pgfxy(80.00,80.00)}\pgflineto{\pgfxy(50.00,80.00)}\pgfclosepath\pgfstroke
\pgfmoveto{\pgfxy(80.00,60.00)}\pgflineto{\pgfxy(95.00,60.00)}\pgflineto{\pgfxy(95.00,65.00)}\pgflineto{\pgfxy(80.00,65.00)}\pgfclosepath\pgfstroke
\pgfmoveto{\pgfxy(130.00,40.00)}\pgflineto{\pgfxy(160.00,40.00)}\pgflineto{\pgfxy(160.00,60.00)}\pgflineto{\pgfxy(130.00,60.00)}\pgfclosepath\pgfstroke
\pgfmoveto{\pgfxy(170.00,40.00)}\pgflineto{\pgfxy(185.00,40.00)}\pgflineto{\pgfxy(185.00,45.00)}\pgflineto{\pgfxy(170.00,45.00)}\pgfclosepath\pgfstroke
\pgfputat{\pgfxy(92.50,61.00)}{\pgfbox[bottom,left]{\fontsize{9.10}{10.93}\selectfont \makebox[0pt]{$E$}}}
\pgfputat{\pgfxy(177.50,41.00)}{\pgfbox[bottom,left]{\fontsize{9.10}{10.93}\selectfont \makebox[0pt]{$E$}}}
\pgfputat{\pgfxy(167.00,97.00)}{\pgfbox[bottom,left]{\fontsize{9.10}{10.93}\selectfont $\tilde{C} = (C,r_t,\ldots,r_1)$}}
\end{pgfpicture}%
$$
\caption{Pointed partition $(\lambda, v)$ and its corresponding quintuples $(\tilde{A},B,\tilde{C},D,E)$ by $\Psi$}
\label{fig:region}
\end{figure}

We can construct the inverse of the mapping $\psi_{a,l,m}$ as follows:
Move out all parts $r_i$, $i=1,\ldots, t$, of size at most $a+m$ from a partition $\tilde{C}$ and denote the remained partition by $C$.
By Lemma~\ref{lem:pealing}, $A$ can be recovered from $(\tilde{A}; r_1, \ldots, r_t)$.

\begin{ex}
Continuing the previous example, we have  $(r_1,r_2,r_3)=(5,7,8)$,
 $(\tilde{A},B,\tilde{C},D,E)\in \mathcal{Q}_{101}(3,3,5)$ with $\tilde{A}=(5,3,1)$, $B=(2,1,1)$, $\tilde{C}=(12,10,10,9,8,7,5)$, $D=(6,6,6,6)$ and $E=(3)$.
\end{ex}

So every pointed partition $(\lambda, v) \in \mathcal{F}_n(a,l,m)$ can be transformed into five partitions $(\tilde{A},B,\tilde{C},D,E) \in \widetilde\mathcal{Q}_n(a,l,m)$ by $\psi_{a,l,m} \circ \varphi_{a,l,m}$, as shown in Figure~\ref{fig:region}.
Let $\widetilde\mathcal{Q}_n$ be the disjoint union of the sets $\widetilde\mathcal{Q}_n(a,l,m)$ for all $a,l,m\geq 0$.
Define the bijection $\Psi$ from $\mathcal{F}_n$ to $\widetilde\mathcal{Q}_n$ by
$$
\Psi(\lambda,v) = \psi_{a,l,m} \circ \varphi_{a,l,m}(\lambda,v) \quad \text{if $(\lambda,v) \in \mathcal{F}_n(a,l,m)$}.
$$

\subsection*{Step 3: The involution $\rho: (\tilde{A},B,\tilde{C},D,E) \mapsto  (B',\tilde{A}',\tilde{C},D',E)$ }
Define the involution $\rho$ on $\widetilde\mathcal{Q}_n$ by
$$
\rho(\tilde{A},B,\tilde{C},D,E) = (B',\tilde{A}',\tilde{C},D',E),
$$
where $X'$ is the conjugate of the partition $X$.

\begin{ex}
Continuing the previous example, we have
$\rho(\tilde{A},B,\tilde{C},D,E) = (B',\tilde{A}',\tilde{C},D',E)$, where
$$
B'=(3,1), \; \tilde{A}'= (3,2,2,1,1),\;  D'=(4,4,4,4,4,4).
$$
\end{ex}
\begin{thm}
\label{thm:main}
For all $n\geq 0$, the mapping $\Phi = \Psi^{-1} \circ \rho \circ \Psi$ is an involution  on $\mathcal{F}_n$ such that
if $\Phi:(\lambda, v) \mapsto (\mu, u)$ then
\begin{equation}
(a_v,l_v,m_v)(\lambda) = (a_u,m_u,l_u)(\mu).
\label{eq:symmetric}
\end{equation}
In particular,  the mapping $\Phi$ also satisfies
\begin{equation}
(h_v,p_v)(\lambda) = (p_u,h_u)(\mu).
\label{eq:symmetric2}
\end{equation}
\end{thm}


\begin{proof}
By definition, the restriction of $\rho$ on $\widetilde\mathcal{Q}_n(a,l,m)$ is  a bijection from $\widetilde\mathcal{Q}_n(a,l,m)$ to $\widetilde\mathcal{Q}_n(a,m,l)$.
Hence, for any  $(\lambda, v) \in \mathcal{F}_n(a,l,m)$, we have
\begin{align*}
\Phi(\lambda,v) &=\Psi^{-1} \circ \rho \circ \Psi(\lambda,v) \\
&= \varphi_{a,m,l}^{-1} \circ \psi_{a,m,l}^{-1} \circ \rho \circ \psi_{a,l,m} \circ \varphi_{a,l,m} (\lambda,v) \in \mathcal{F}_n(a,m,l).
\end{align*}
Clearly the mapping $\Phi$ is an involution on $\mathcal{F}_n$ satisfying \eqref{eq:symmetric} and \eqref{eq:symmetric2}.
\end{proof}

\rmk In other words,  we have the following diagram:
$$
\vcenter{
\xymatrix{
&\mathcal{F}_n(a,l,m) \ar[rr]^-{\Phi} \ar[ld]_{\varphi_{a,l,m}} \ar[dd]_{\Psi} &&\mathcal{F}_n(a,m,l)& \\
\mathcal{Q}_n(a,l,m) \ar[rd]_{\psi_{a,l,m}}& & &&\mathcal{Q}_n(a,m,l) \ar[lu]_{\varphi_{a,m,l}^{-1}}\\
&\widetilde\mathcal{Q}_n(a,l,m) \ar[rr]^-{\rho} &&\widetilde\mathcal{Q}_n(a,m,l) \ar[uu]_{\Psi^{-1}} \ar[ru]_{\psi_{a,m,l}^{-1}}
}}
.
$$

\begin{ex}
Continuing the above example, we have $\Phi(\lambda, v)=(\mu,u) \in \mathcal{F}_{101}(3,5,3)$ with
$$\mu=(12,10,10,9,8,7,7,7,6,6,5,5,3,2,2,1,1)\quad \text{and}\quad  u=(4,7).$$
\end{ex}

To end this section we draw a graph on ${\mathcal F}_4$ to illustrate the bijection $\Phi$ on $\mathcal{F}_4$ by connecting each $(\lambda, v)$ to $\Phi(\lambda,v)$ (see Figure~\ref{fig:part}) as follows:
$$
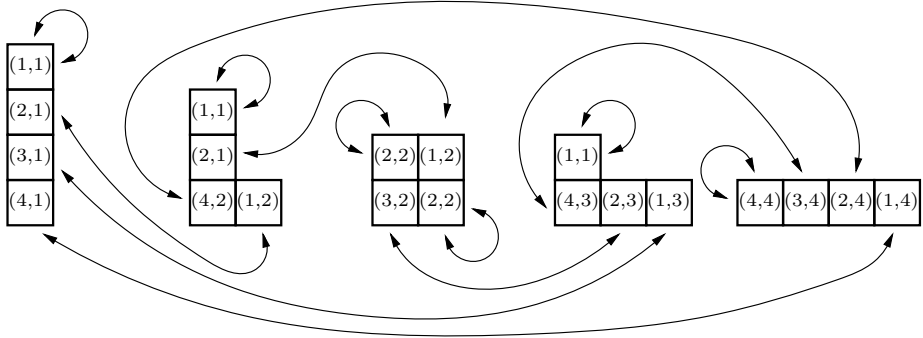
\begin{pgfpicture}{4.00mm}{19.18mm}{128.00mm}{67.53mm}
\pgfsetxvec{\pgfpoint{0.60mm}{0mm}}
\pgfsetyvec{\pgfpoint{0mm}{0.60mm}}
\color[rgb]{0,0,0}\pgfsetlinewidth{0.30mm}\pgfsetdash{}{0mm}
\pgfmoveto{\pgfxy(10.00,89.74)}\pgflineto{\pgfxy(20.00,89.74)}\pgflineto{\pgfxy(20.00,99.74)}\pgflineto{\pgfxy(10.00,99.74)}\pgfclosepath\pgfstroke
\pgfmoveto{\pgfxy(50.00,79.74)}\pgflineto{\pgfxy(60.00,79.74)}\pgflineto{\pgfxy(60.00,89.74)}\pgflineto{\pgfxy(50.00,89.74)}\pgfclosepath\pgfstroke
\pgfmoveto{\pgfxy(10.00,79.74)}\pgflineto{\pgfxy(20.00,79.74)}\pgflineto{\pgfxy(20.00,89.74)}\pgflineto{\pgfxy(10.00,89.74)}\pgfclosepath\pgfstroke
\pgfmoveto{\pgfxy(10.00,69.74)}\pgflineto{\pgfxy(20.00,69.74)}\pgflineto{\pgfxy(20.00,79.74)}\pgflineto{\pgfxy(10.00,79.74)}\pgfclosepath\pgfstroke
\pgfmoveto{\pgfxy(10.00,59.74)}\pgflineto{\pgfxy(20.00,59.74)}\pgflineto{\pgfxy(20.00,69.74)}\pgflineto{\pgfxy(10.00,69.74)}\pgfclosepath\pgfstroke
\pgfmoveto{\pgfxy(50.00,69.74)}\pgflineto{\pgfxy(60.00,69.74)}\pgflineto{\pgfxy(60.00,79.74)}\pgflineto{\pgfxy(50.00,79.74)}\pgfclosepath\pgfstroke
\pgfmoveto{\pgfxy(50.00,59.74)}\pgflineto{\pgfxy(60.00,59.74)}\pgflineto{\pgfxy(60.00,69.74)}\pgflineto{\pgfxy(50.00,69.74)}\pgfclosepath\pgfstroke
\pgfmoveto{\pgfxy(60.00,59.74)}\pgflineto{\pgfxy(70.00,59.74)}\pgflineto{\pgfxy(70.00,69.74)}\pgflineto{\pgfxy(60.00,69.74)}\pgfclosepath\pgfstroke
\pgfmoveto{\pgfxy(90.00,69.74)}\pgflineto{\pgfxy(100.00,69.74)}\pgflineto{\pgfxy(100.00,79.74)}\pgflineto{\pgfxy(90.00,79.74)}\pgfclosepath\pgfstroke
\pgfmoveto{\pgfxy(90.00,59.74)}\pgflineto{\pgfxy(100.00,59.74)}\pgflineto{\pgfxy(100.00,69.74)}\pgflineto{\pgfxy(90.00,69.74)}\pgfclosepath\pgfstroke
\pgfmoveto{\pgfxy(100.00,69.74)}\pgflineto{\pgfxy(110.00,69.74)}\pgflineto{\pgfxy(110.00,79.74)}\pgflineto{\pgfxy(100.00,79.74)}\pgfclosepath\pgfstroke
\pgfmoveto{\pgfxy(100.00,59.74)}\pgflineto{\pgfxy(110.00,59.74)}\pgflineto{\pgfxy(110.00,69.74)}\pgflineto{\pgfxy(100.00,69.74)}\pgfclosepath\pgfstroke
\pgfmoveto{\pgfxy(130.00,69.74)}\pgflineto{\pgfxy(140.00,69.74)}\pgflineto{\pgfxy(140.00,79.74)}\pgflineto{\pgfxy(130.00,79.74)}\pgfclosepath\pgfstroke
\pgfmoveto{\pgfxy(130.00,59.74)}\pgflineto{\pgfxy(140.00,59.74)}\pgflineto{\pgfxy(140.00,69.74)}\pgflineto{\pgfxy(130.00,69.74)}\pgfclosepath\pgfstroke
\pgfmoveto{\pgfxy(140.00,59.74)}\pgflineto{\pgfxy(150.00,59.74)}\pgflineto{\pgfxy(150.00,69.74)}\pgflineto{\pgfxy(140.00,69.74)}\pgfclosepath\pgfstroke
\pgfmoveto{\pgfxy(150.00,59.74)}\pgflineto{\pgfxy(160.00,59.74)}\pgflineto{\pgfxy(160.00,69.74)}\pgflineto{\pgfxy(150.00,69.74)}\pgfclosepath\pgfstroke
\pgfmoveto{\pgfxy(170.00,59.74)}\pgflineto{\pgfxy(180.00,59.74)}\pgflineto{\pgfxy(180.00,69.74)}\pgflineto{\pgfxy(170.00,69.74)}\pgfclosepath\pgfstroke
\pgfmoveto{\pgfxy(180.00,59.74)}\pgflineto{\pgfxy(190.00,59.74)}\pgflineto{\pgfxy(190.00,69.74)}\pgflineto{\pgfxy(180.00,69.74)}\pgfclosepath\pgfstroke
\pgfmoveto{\pgfxy(190.00,59.74)}\pgflineto{\pgfxy(200.00,59.74)}\pgflineto{\pgfxy(200.00,69.74)}\pgflineto{\pgfxy(190.00,69.74)}\pgfclosepath\pgfstroke
\pgfmoveto{\pgfxy(200.00,59.74)}\pgflineto{\pgfxy(210.00,59.74)}\pgflineto{\pgfxy(210.00,69.74)}\pgflineto{\pgfxy(200.00,69.74)}\pgfclosepath\pgfstroke
\pgfsetlinewidth{0.15mm}\pgfmoveto{\pgfxy(18.00,57.74)}\pgfcurveto{\pgfxy(31.59,49.39)}{\pgfxy(46.45,43.31)}{\pgfxy(62.00,39.74)}\pgfcurveto{\pgfxy(82.94,34.95)}{\pgfxy(104.55,34.82)}{\pgfxy(126.00,35.74)}\pgfcurveto{\pgfxy(138.74,36.29)}{\pgfxy(151.49,37.22)}{\pgfxy(164.00,39.74)}\pgfcurveto{\pgfxy(173.53,41.67)}{\pgfxy(182.83,44.52)}{\pgfxy(192.00,47.74)}\pgfcurveto{\pgfxy(194.84,48.74)}{\pgfxy(197.70,49.80)}{\pgfxy(200.00,51.74)}\pgfcurveto{\pgfxy(201.87,53.33)}{\pgfxy(203.26,55.41)}{\pgfxy(204.00,57.74)}\pgfstroke
\pgfmoveto{\pgfxy(18.00,57.74)}\pgflineto{\pgfxy(20.04,55.71)}\pgflineto{\pgfxy(20.76,56.91)}\pgflineto{\pgfxy(18.00,57.74)}\pgfclosepath\pgffill
\pgfmoveto{\pgfxy(18.00,57.74)}\pgflineto{\pgfxy(20.04,55.71)}\pgflineto{\pgfxy(20.76,56.91)}\pgflineto{\pgfxy(18.00,57.74)}\pgfclosepath\pgfstroke
\pgfmoveto{\pgfxy(204.00,57.74)}\pgflineto{\pgfxy(202.15,55.53)}\pgflineto{\pgfxy(203.42,54.92)}\pgflineto{\pgfxy(204.00,57.74)}\pgfclosepath\pgffill
\pgfmoveto{\pgfxy(204.00,57.74)}\pgflineto{\pgfxy(202.15,55.53)}\pgflineto{\pgfxy(203.42,54.92)}\pgflineto{\pgfxy(204.00,57.74)}\pgfclosepath\pgfstroke
\pgfmoveto{\pgfxy(22.00,71.74)}\pgfcurveto{\pgfxy(32.24,59.98)}{\pgfxy(45.32,51.03)}{\pgfxy(60.00,45.74)}\pgfcurveto{\pgfxy(69.02,42.50)}{\pgfxy(78.47,40.70)}{\pgfxy(88.00,39.74)}\pgfcurveto{\pgfxy(100.05,38.53)}{\pgfxy(112.28,38.66)}{\pgfxy(124.00,41.74)}\pgfcurveto{\pgfxy(135.11,44.66)}{\pgfxy(145.39,50.14)}{\pgfxy(154.00,57.74)}\pgfstroke
\pgfmoveto{\pgfxy(22.00,71.74)}\pgflineto{\pgfxy(23.36,69.20)}\pgflineto{\pgfxy(24.40,70.14)}\pgflineto{\pgfxy(22.00,71.74)}\pgfclosepath\pgffill
\pgfmoveto{\pgfxy(22.00,71.74)}\pgflineto{\pgfxy(23.36,69.20)}\pgflineto{\pgfxy(24.40,70.14)}\pgflineto{\pgfxy(22.00,71.74)}\pgfclosepath\pgfstroke
\pgfmoveto{\pgfxy(154.00,57.74)}\pgflineto{\pgfxy(151.40,56.48)}\pgflineto{\pgfxy(152.30,55.41)}\pgflineto{\pgfxy(154.00,57.74)}\pgfclosepath\pgffill
\pgfmoveto{\pgfxy(154.00,57.74)}\pgflineto{\pgfxy(151.40,56.48)}\pgflineto{\pgfxy(152.30,55.41)}\pgflineto{\pgfxy(154.00,57.74)}\pgfclosepath\pgfstroke
\pgfmoveto{\pgfxy(22.00,83.74)}\pgfcurveto{\pgfxy(26.29,76.30)}{\pgfxy(31.68,69.56)}{\pgfxy(38.00,63.74)}\pgfcurveto{\pgfxy(42.92,59.22)}{\pgfxy(48.34,55.30)}{\pgfxy(54.00,51.74)}\pgfcurveto{\pgfxy(57.15,49.77)}{\pgfxy(60.81,48.02)}{\pgfxy(64.00,49.74)}\pgfcurveto{\pgfxy(66.88,51.31)}{\pgfxy(67.81,55.01)}{\pgfxy(66.00,57.74)}\pgfstroke
\pgfmoveto{\pgfxy(22.00,83.74)}\pgflineto{\pgfxy(22.85,80.99)}\pgflineto{\pgfxy(24.05,81.71)}\pgflineto{\pgfxy(22.00,83.74)}\pgfclosepath\pgffill
\pgfmoveto{\pgfxy(22.00,83.74)}\pgflineto{\pgfxy(22.85,80.99)}\pgflineto{\pgfxy(24.05,81.71)}\pgflineto{\pgfxy(22.00,83.74)}\pgfclosepath\pgfstroke
\pgfmoveto{\pgfxy(66.00,57.74)}\pgflineto{\pgfxy(66.03,54.86)}\pgflineto{\pgfxy(67.38,55.21)}\pgflineto{\pgfxy(66.00,57.74)}\pgfclosepath\pgffill
\pgfmoveto{\pgfxy(66.00,57.74)}\pgflineto{\pgfxy(66.03,54.86)}\pgflineto{\pgfxy(67.38,55.21)}\pgflineto{\pgfxy(66.00,57.74)}\pgfclosepath\pgfstroke
\pgfmoveto{\pgfxy(16.00,101.74)}\pgfcurveto{\pgfxy(16.13,103.28)}{\pgfxy(16.85,104.71)}{\pgfxy(18.00,105.74)}\pgfcurveto{\pgfxy(20.32,107.82)}{\pgfxy(23.83,107.91)}{\pgfxy(26.00,105.74)}\pgfcurveto{\pgfxy(28.17,103.58)}{\pgfxy(28.08,100.07)}{\pgfxy(26.00,97.74)}\pgfcurveto{\pgfxy(24.97,96.59)}{\pgfxy(23.54,95.88)}{\pgfxy(22.00,95.74)}\pgfstroke
\pgfmoveto{\pgfxy(16.00,101.74)}\pgflineto{\pgfxy(17.72,104.06)}\pgflineto{\pgfxy(16.43,104.60)}\pgflineto{\pgfxy(16.00,101.74)}\pgfclosepath\pgffill
\pgfmoveto{\pgfxy(16.00,101.74)}\pgflineto{\pgfxy(17.72,104.06)}\pgflineto{\pgfxy(16.43,104.60)}\pgflineto{\pgfxy(16.00,101.74)}\pgfclosepath\pgfstroke
\pgfmoveto{\pgfxy(22.00,95.74)}\pgflineto{\pgfxy(24.85,96.17)}\pgflineto{\pgfxy(24.32,97.46)}\pgflineto{\pgfxy(22.00,95.74)}\pgfclosepath\pgffill
\pgfmoveto{\pgfxy(22.00,95.74)}\pgflineto{\pgfxy(24.85,96.17)}\pgflineto{\pgfxy(24.32,97.46)}\pgflineto{\pgfxy(22.00,95.74)}\pgfclosepath\pgfstroke
\pgfmoveto{\pgfxy(56.00,91.74)}\pgfcurveto{\pgfxy(56.13,93.28)}{\pgfxy(56.85,94.71)}{\pgfxy(58.00,95.74)}\pgfcurveto{\pgfxy(60.32,97.82)}{\pgfxy(63.83,97.91)}{\pgfxy(66.00,95.74)}\pgfcurveto{\pgfxy(68.17,93.58)}{\pgfxy(68.08,90.07)}{\pgfxy(66.00,87.74)}\pgfcurveto{\pgfxy(64.97,86.59)}{\pgfxy(63.54,85.88)}{\pgfxy(62.00,85.74)}\pgfstroke
\pgfmoveto{\pgfxy(56.00,91.74)}\pgflineto{\pgfxy(57.72,94.06)}\pgflineto{\pgfxy(56.43,94.60)}\pgflineto{\pgfxy(56.00,91.74)}\pgfclosepath\pgffill
\pgfmoveto{\pgfxy(56.00,91.74)}\pgflineto{\pgfxy(57.72,94.06)}\pgflineto{\pgfxy(56.43,94.60)}\pgflineto{\pgfxy(56.00,91.74)}\pgfclosepath\pgfstroke
\pgfmoveto{\pgfxy(62.00,85.74)}\pgflineto{\pgfxy(64.85,86.17)}\pgflineto{\pgfxy(64.32,87.46)}\pgflineto{\pgfxy(62.00,85.74)}\pgfclosepath\pgffill
\pgfmoveto{\pgfxy(62.00,85.74)}\pgflineto{\pgfxy(64.85,86.17)}\pgflineto{\pgfxy(64.32,87.46)}\pgflineto{\pgfxy(62.00,85.74)}\pgfclosepath\pgfstroke
\pgfmoveto{\pgfxy(136.00,81.74)}\pgfcurveto{\pgfxy(136.13,83.28)}{\pgfxy(136.85,84.71)}{\pgfxy(138.00,85.74)}\pgfcurveto{\pgfxy(140.32,87.82)}{\pgfxy(143.83,87.91)}{\pgfxy(146.00,85.74)}\pgfcurveto{\pgfxy(148.17,83.58)}{\pgfxy(148.08,80.07)}{\pgfxy(146.00,77.74)}\pgfcurveto{\pgfxy(144.97,76.59)}{\pgfxy(143.54,75.88)}{\pgfxy(142.00,75.74)}\pgfstroke
\pgfmoveto{\pgfxy(136.00,81.74)}\pgflineto{\pgfxy(137.72,84.06)}\pgflineto{\pgfxy(136.43,84.60)}\pgflineto{\pgfxy(136.00,81.74)}\pgfclosepath\pgffill
\pgfmoveto{\pgfxy(136.00,81.74)}\pgflineto{\pgfxy(137.72,84.06)}\pgflineto{\pgfxy(136.43,84.60)}\pgflineto{\pgfxy(136.00,81.74)}\pgfclosepath\pgfstroke
\pgfmoveto{\pgfxy(142.00,75.74)}\pgflineto{\pgfxy(144.85,76.17)}\pgflineto{\pgfxy(144.32,77.46)}\pgflineto{\pgfxy(142.00,75.74)}\pgfclosepath\pgffill
\pgfmoveto{\pgfxy(142.00,75.74)}\pgflineto{\pgfxy(144.85,76.17)}\pgflineto{\pgfxy(144.32,77.46)}\pgflineto{\pgfxy(142.00,75.74)}\pgfclosepath\pgfstroke
\pgfmoveto{\pgfxy(93.59,81.74)}\pgfcurveto{\pgfxy(93.46,83.28)}{\pgfxy(92.74,84.71)}{\pgfxy(91.59,85.74)}\pgfcurveto{\pgfxy(89.27,87.82)}{\pgfxy(85.76,87.91)}{\pgfxy(83.59,85.74)}\pgfcurveto{\pgfxy(81.43,83.58)}{\pgfxy(81.51,80.07)}{\pgfxy(83.59,77.74)}\pgfcurveto{\pgfxy(84.62,76.59)}{\pgfxy(86.05,75.88)}{\pgfxy(87.59,75.74)}\pgfstroke
\pgfmoveto{\pgfxy(93.59,81.74)}\pgflineto{\pgfxy(93.17,84.60)}\pgflineto{\pgfxy(91.87,84.06)}\pgflineto{\pgfxy(93.59,81.74)}\pgfclosepath\pgffill
\pgfmoveto{\pgfxy(93.59,81.74)}\pgflineto{\pgfxy(93.17,84.60)}\pgflineto{\pgfxy(91.87,84.06)}\pgflineto{\pgfxy(93.59,81.74)}\pgfclosepath\pgfstroke
\pgfmoveto{\pgfxy(87.59,75.74)}\pgflineto{\pgfxy(85.27,77.46)}\pgflineto{\pgfxy(84.74,76.17)}\pgflineto{\pgfxy(87.59,75.74)}\pgfclosepath\pgffill
\pgfmoveto{\pgfxy(87.59,75.74)}\pgflineto{\pgfxy(85.27,77.46)}\pgflineto{\pgfxy(84.74,76.17)}\pgflineto{\pgfxy(87.59,75.74)}\pgfclosepath\pgfstroke
\pgfmoveto{\pgfxy(106.00,57.34)}\pgfcurveto{\pgfxy(106.13,55.80)}{\pgfxy(106.85,54.37)}{\pgfxy(108.00,53.34)}\pgfcurveto{\pgfxy(110.32,51.26)}{\pgfxy(113.83,51.17)}{\pgfxy(116.00,53.34)}\pgfcurveto{\pgfxy(118.17,55.50)}{\pgfxy(118.08,59.01)}{\pgfxy(116.00,61.34)}\pgfcurveto{\pgfxy(114.97,62.49)}{\pgfxy(113.54,63.20)}{\pgfxy(112.00,63.34)}\pgfstroke
\pgfmoveto{\pgfxy(106.00,57.34)}\pgflineto{\pgfxy(106.43,54.48)}\pgflineto{\pgfxy(107.72,55.02)}\pgflineto{\pgfxy(106.00,57.34)}\pgfclosepath\pgffill
\pgfmoveto{\pgfxy(106.00,57.34)}\pgflineto{\pgfxy(106.43,54.48)}\pgflineto{\pgfxy(107.72,55.02)}\pgflineto{\pgfxy(106.00,57.34)}\pgfclosepath\pgfstroke
\pgfmoveto{\pgfxy(112.00,63.34)}\pgflineto{\pgfxy(114.32,61.62)}\pgflineto{\pgfxy(114.85,62.91)}\pgflineto{\pgfxy(112.00,63.34)}\pgfclosepath\pgffill
\pgfmoveto{\pgfxy(112.00,63.34)}\pgflineto{\pgfxy(114.32,61.62)}\pgflineto{\pgfxy(114.85,62.91)}\pgflineto{\pgfxy(112.00,63.34)}\pgfclosepath\pgfstroke
\pgfmoveto{\pgfxy(94.00,57.74)}\pgfcurveto{\pgfxy(94.91,55.49)}{\pgfxy(96.27,53.45)}{\pgfxy(98.00,51.74)}\pgfcurveto{\pgfxy(101.73,48.08)}{\pgfxy(106.79,46.23)}{\pgfxy(112.00,45.74)}\pgfcurveto{\pgfxy(118.90,45.10)}{\pgfxy(125.75,46.77)}{\pgfxy(132.00,49.74)}\pgfcurveto{\pgfxy(136.36,51.82)}{\pgfxy(140.40,54.52)}{\pgfxy(144.00,57.74)}\pgfstroke
\pgfmoveto{\pgfxy(94.00,57.74)}\pgflineto{\pgfxy(94.69,54.94)}\pgflineto{\pgfxy(95.93,55.60)}\pgflineto{\pgfxy(94.00,57.74)}\pgfclosepath\pgffill
\pgfmoveto{\pgfxy(94.00,57.74)}\pgflineto{\pgfxy(94.69,54.94)}\pgflineto{\pgfxy(95.93,55.60)}\pgflineto{\pgfxy(94.00,57.74)}\pgfclosepath\pgfstroke
\pgfmoveto{\pgfxy(144.00,57.74)}\pgflineto{\pgfxy(141.40,56.49)}\pgflineto{\pgfxy(142.30,55.41)}\pgflineto{\pgfxy(144.00,57.74)}\pgfclosepath\pgffill
\pgfmoveto{\pgfxy(144.00,57.74)}\pgflineto{\pgfxy(141.40,56.49)}\pgflineto{\pgfxy(142.30,55.41)}\pgflineto{\pgfxy(144.00,57.74)}\pgfclosepath\pgfstroke
\pgfmoveto{\pgfxy(128.00,63.74)}\pgfcurveto{\pgfxy(122.10,69.47)}{\pgfxy(120.50,78.31)}{\pgfxy(124.00,85.74)}\pgfcurveto{\pgfxy(126.71,91.51)}{\pgfxy(132.07,95.40)}{\pgfxy(138.00,97.74)}\pgfcurveto{\pgfxy(142.46,99.51)}{\pgfxy(147.25,100.42)}{\pgfxy(152.00,99.74)}\pgfcurveto{\pgfxy(158.97,98.76)}{\pgfxy(164.89,94.58)}{\pgfxy(170.00,89.74)}\pgfcurveto{\pgfxy(175.55,84.49)}{\pgfxy(180.27,78.42)}{\pgfxy(184.00,71.74)}\pgfstroke
\pgfmoveto{\pgfxy(128.00,63.74)}\pgflineto{\pgfxy(126.72,66.33)}\pgflineto{\pgfxy(125.65,65.42)}\pgflineto{\pgfxy(128.00,63.74)}\pgfclosepath\pgffill
\pgfmoveto{\pgfxy(128.00,63.74)}\pgflineto{\pgfxy(126.72,66.33)}\pgflineto{\pgfxy(125.65,65.42)}\pgflineto{\pgfxy(128.00,63.74)}\pgfclosepath\pgfstroke
\pgfmoveto{\pgfxy(184.00,71.74)}\pgflineto{\pgfxy(183.18,74.51)}\pgflineto{\pgfxy(181.97,73.80)}\pgflineto{\pgfxy(184.00,71.74)}\pgfclosepath\pgffill
\pgfmoveto{\pgfxy(184.00,71.74)}\pgflineto{\pgfxy(183.18,74.51)}\pgflineto{\pgfxy(181.97,73.80)}\pgflineto{\pgfxy(184.00,71.74)}\pgfclosepath\pgfstroke
\pgfmoveto{\pgfxy(48.00,65.74)}\pgfcurveto{\pgfxy(37.58,68.50)}{\pgfxy(32.62,80.40)}{\pgfxy(38.00,89.74)}\pgfcurveto{\pgfxy(40.95,94.87)}{\pgfxy(46.49,97.60)}{\pgfxy(52.00,99.74)}\pgfcurveto{\pgfxy(88.25,113.84)}{\pgfxy(128.07,110.22)}{\pgfxy(166.00,101.74)}\pgfcurveto{\pgfxy(178.84,98.87)}{\pgfxy(192.39,94.06)}{\pgfxy(196.00,81.74)}\pgfcurveto{\pgfxy(196.96,78.48)}{\pgfxy(196.96,75.01)}{\pgfxy(196.00,71.74)}\pgfstroke
\pgfmoveto{\pgfxy(48.00,65.74)}\pgflineto{\pgfxy(45.64,67.41)}\pgflineto{\pgfxy(45.14,66.10)}\pgflineto{\pgfxy(48.00,65.74)}\pgfclosepath\pgffill
\pgfmoveto{\pgfxy(48.00,65.74)}\pgflineto{\pgfxy(45.64,67.41)}\pgflineto{\pgfxy(45.14,66.10)}\pgflineto{\pgfxy(48.00,65.74)}\pgfclosepath\pgfstroke
\pgfmoveto{\pgfxy(196.00,71.74)}\pgflineto{\pgfxy(197.19,74.37)}\pgflineto{\pgfxy(195.82,74.62)}\pgflineto{\pgfxy(196.00,71.74)}\pgfclosepath\pgffill
\pgfmoveto{\pgfxy(196.00,71.74)}\pgflineto{\pgfxy(197.19,74.37)}\pgflineto{\pgfxy(195.82,74.62)}\pgflineto{\pgfxy(196.00,71.74)}\pgfclosepath\pgfstroke
\pgfmoveto{\pgfxy(173.59,71.74)}\pgfcurveto{\pgfxy(173.46,73.28)}{\pgfxy(172.74,74.71)}{\pgfxy(171.59,75.74)}\pgfcurveto{\pgfxy(169.27,77.82)}{\pgfxy(165.76,77.91)}{\pgfxy(163.59,75.74)}\pgfcurveto{\pgfxy(161.43,73.58)}{\pgfxy(161.51,70.07)}{\pgfxy(163.59,67.74)}\pgfcurveto{\pgfxy(164.62,66.59)}{\pgfxy(166.05,65.88)}{\pgfxy(167.59,65.74)}\pgfstroke
\pgfmoveto{\pgfxy(173.59,71.74)}\pgflineto{\pgfxy(173.17,74.60)}\pgflineto{\pgfxy(171.87,74.06)}\pgflineto{\pgfxy(173.59,71.74)}\pgfclosepath\pgffill
\pgfmoveto{\pgfxy(173.59,71.74)}\pgflineto{\pgfxy(173.17,74.60)}\pgflineto{\pgfxy(171.87,74.06)}\pgflineto{\pgfxy(173.59,71.74)}\pgfclosepath\pgfstroke
\pgfmoveto{\pgfxy(167.59,65.74)}\pgflineto{\pgfxy(165.27,67.46)}\pgflineto{\pgfxy(164.74,66.17)}\pgflineto{\pgfxy(167.59,65.74)}\pgfclosepath\pgffill
\pgfmoveto{\pgfxy(167.59,65.74)}\pgflineto{\pgfxy(165.27,67.46)}\pgflineto{\pgfxy(164.74,66.17)}\pgflineto{\pgfxy(167.59,65.74)}\pgfclosepath\pgfstroke
\pgfmoveto{\pgfxy(62.00,75.74)}\pgfcurveto{\pgfxy(63.97,75.24)}{\pgfxy(66.03,75.24)}{\pgfxy(68.00,75.74)}\pgfcurveto{\pgfxy(71.35,76.60)}{\pgfxy(74.17,78.81)}{\pgfxy(76.00,81.74)}\pgfcurveto{\pgfxy(78.39,85.57)}{\pgfxy(78.96,90.39)}{\pgfxy(82.00,93.74)}\pgfcurveto{\pgfxy(86.46,98.67)}{\pgfxy(93.92,98.86)}{\pgfxy(100.00,95.74)}\pgfcurveto{\pgfxy(102.61,94.41)}{\pgfxy(104.95,92.47)}{\pgfxy(106.00,89.74)}\pgfcurveto{\pgfxy(106.75,87.81)}{\pgfxy(106.75,85.67)}{\pgfxy(106.00,83.74)}\pgfstroke
\pgfmoveto{\pgfxy(62.00,75.74)}\pgflineto{\pgfxy(64.72,74.78)}\pgflineto{\pgfxy(64.85,76.17)}\pgflineto{\pgfxy(62.00,75.74)}\pgfclosepath\pgffill
\pgfmoveto{\pgfxy(62.00,75.74)}\pgflineto{\pgfxy(64.72,74.78)}\pgflineto{\pgfxy(64.85,76.17)}\pgflineto{\pgfxy(62.00,75.74)}\pgfclosepath\pgfstroke
\pgfmoveto{\pgfxy(106.00,83.74)}\pgflineto{\pgfxy(107.10,86.41)}\pgflineto{\pgfxy(105.72,86.62)}\pgflineto{\pgfxy(106.00,83.74)}\pgfclosepath\pgffill
\pgfmoveto{\pgfxy(106.00,83.74)}\pgflineto{\pgfxy(107.10,86.41)}\pgflineto{\pgfxy(105.72,86.62)}\pgflineto{\pgfxy(106.00,83.74)}\pgfclosepath\pgfstroke
\pgfputat{\pgfxy(15.00,94.00)}{\pgfbox[bottom,left]{\fontsize{6.83}{8.19}\selectfont \makebox[0pt]{(1,1)}}}
\pgfputat{\pgfxy(55.00,84.00)}{\pgfbox[bottom,left]{\fontsize{6.83}{8.19}\selectfont \makebox[0pt]{(1,1)}}}
\pgfputat{\pgfxy(55.00,74.00)}{\pgfbox[bottom,left]{\fontsize{6.83}{8.19}\selectfont \makebox[0pt]{(2,1)}}}
\pgfputat{\pgfxy(55.00,64.00)}{\pgfbox[bottom,left]{\fontsize{6.83}{8.19}\selectfont \makebox[0pt]{(4,2)}}}
\pgfputat{\pgfxy(15.00,84.00)}{\pgfbox[bottom,left]{\fontsize{6.83}{8.19}\selectfont \makebox[0pt]{(2,1)}}}
\pgfputat{\pgfxy(15.00,74.00)}{\pgfbox[bottom,left]{\fontsize{6.83}{8.19}\selectfont \makebox[0pt]{(3,1)}}}
\pgfputat{\pgfxy(15.00,64.00)}{\pgfbox[bottom,left]{\fontsize{6.83}{8.19}\selectfont \makebox[0pt]{(4,1)}}}
\pgfputat{\pgfxy(65.00,64.00)}{\pgfbox[bottom,left]{\fontsize{6.83}{8.19}\selectfont \makebox[0pt]{(1,2)}}}
\pgfputat{\pgfxy(95.00,74.00)}{\pgfbox[bottom,left]{\fontsize{6.83}{8.19}\selectfont \makebox[0pt]{(2,2)}}}
\pgfputat{\pgfxy(105.00,74.00)}{\pgfbox[bottom,left]{\fontsize{6.83}{8.19}\selectfont \makebox[0pt]{(1,2)}}}
\pgfputat{\pgfxy(95.00,64.00)}{\pgfbox[bottom,left]{\fontsize{6.83}{8.19}\selectfont \makebox[0pt]{(3,2)}}}
\pgfputat{\pgfxy(105.00,64.00)}{\pgfbox[bottom,left]{\fontsize{6.83}{8.19}\selectfont \makebox[0pt]{(2,2)}}}
\pgfputat{\pgfxy(135.00,74.00)}{\pgfbox[bottom,left]{\fontsize{6.83}{8.19}\selectfont \makebox[0pt]{(1,1)}}}
\pgfputat{\pgfxy(155.00,64.00)}{\pgfbox[bottom,left]{\fontsize{6.83}{8.19}\selectfont \makebox[0pt]{(1,3)}}}
\pgfputat{\pgfxy(145.00,64.00)}{\pgfbox[bottom,left]{\fontsize{6.83}{8.19}\selectfont \makebox[0pt]{(2,3)}}}
\pgfputat{\pgfxy(135.00,64.00)}{\pgfbox[bottom,left]{\fontsize{6.83}{8.19}\selectfont \makebox[0pt]{(4,3)}}}
\pgfputat{\pgfxy(175.00,64.00)}{\pgfbox[bottom,left]{\fontsize{6.83}{8.19}\selectfont \makebox[0pt]{(4,4)}}}
\pgfputat{\pgfxy(185.00,64.00)}{\pgfbox[bottom,left]{\fontsize{6.83}{8.19}\selectfont \makebox[0pt]{(3,4)}}}
\pgfputat{\pgfxy(195.00,64.00)}{\pgfbox[bottom,left]{\fontsize{6.83}{8.19}\selectfont \makebox[0pt]{(2,4)}}}
\pgfputat{\pgfxy(205.00,64.00)}{\pgfbox[bottom,left]{\fontsize{6.83}{8.19}\selectfont \makebox[0pt]{(1,4)}}}
\end{pgfpicture}%
.$$

\section{Some consequences}
We derive immediately the following  result of Bessenrodt and Han \cite[Theorem 3]{BESSENRODT:2009:HAL-00395700:1}.
\begin{cor} The triple statistic $(a_v,l_v,m_v)$ has the same distribution as $(a_v, m_v,l_v)$ over ${\mathcal F}_n$. In other words, the polynomial
$$
Q_n(x,y,z)=\sum_{(\lambda,v)\in {\cal F}_n}x^{a_v}y^{l_v}z^{m_v}
$$
is symmetric in $y$ an $z$.
\end{cor}

Let $f_n(a,l,m)$ be the cardinality of ${\cal F}_n(a,l,m)$, that is, $f_n(a,l,m)$ is the coefficient of $x^ay^lz^m$ in $Q_n(x,y,z)$.
We can apply the bijection $\Phi$ to give a different proof of Bessenrodt and Han's formula \cite[Theorem 2]{BESSENRODT:2009:HAL-00395700:1} for $\sum_{n\geq 0}f_n(a,l,m) q^n$.
\begin{cor}
The generating function of $f_n(a,l,m)$ is given by the following formula
\begin{equation*}\label{eq:gf}
\sum_{n\geq 0}f_n(a,l,m) q^n = \frac{1}{(q^{a+1};q)_{\infty}} \qbin{m+a}{a} \qbin{l+a}{a} q^{(m+1)(l+1)+a},
\end{equation*}
where $(a;q)_{\infty} := \prod_{i=0}^{\infty} (1-aq^i)$.
\end{cor}
\begin{proof} By the bijection $\psi$
the generating function
$\sum_{\lambda} q^{\abs{\lambda}}$, where $(\lambda,v)\in \mathcal{F}_n(a,l,m)$,  is equal to the product of the corresponding
generating functions of the  five partitions $A$, $B$, $\tilde{C}$, $D$ and $E$. In view of the basic facts about partitions (see \cite[Chapter~3]{MR1634067}), it is easy to see that $D(q) = q^{(m+1)(l+1)}$,
$E(q)= q^{a}$, and
\begin{align*}
A(q) = \qbin{m+a}{a};\quad
B(q) = \qbin{l+a}{a};\quad
\tilde{C}(q) = \frac{1}{(q^{a+1};q)_{\infty}}.
\end{align*}
Multiplying the five generating functions yields the result.
\end{proof}

A polynomial $P(x,y)$ in two variables $x$ and $y$ is {\em super-symmetric} if
$$[x^{\alpha} y^{\beta}]P(x,y) = [x^{\alpha'} y^{\beta'}]P(x,y)$$
when $\alpha+\beta=\alpha'+\beta'$. Clearly,  any super-symmetric polynomial  is also symmetric.

It is known (see \cite{MR1682922,MR1894024, BESSENRODT:2009:HAL-00395700:1}) that
the generating function for the pointed partitions of  $\mathcal{F}_n$ by the  two joint statistics arm length and coarm length (resp. leg length)
is super-symmetric.  In other words, the polynomial
$$
\sum_{(\lambda,v)\in \mathcal{F}_n} x^{a_v} y^{m_v}\quad (\textrm{resp.}\quad \sum_{(\lambda,v)\in \mathcal{F}_n} x^{a_v} y^{l_v})
$$
is super-symmetric. Note that the above two polynomials are actually equal due to Corollary~2.

Let $\mathcal{F}_n(a,*,m)$ (resp. $\mathcal{F}_n(a,l,*)$) be the set of pointed partitions $(\lambda,v)$ of $n$ such that $a_v=a$ and $m_v=m$ (resp. $a_v=a$ and $l_v=l$).

It is easy to give a combinatorial proof of the super-symmetry of the first polynomial $\sum_{(\lambda,v)\in \mathcal{F}_n} x^{a_v} y^{m_v}$. Indeed, if  $\alpha+\beta=\alpha'+\beta'$ and $(\lambda, v)\in \mathcal{F}_n(\alpha,*,\beta)$,
let  $u$ be the unique cell on  the same row as $v$ satisfying $a_u(\lambda)=\alpha'$ and $m_u(\lambda)=\beta'$,
then $\tau_{\alpha,\beta,\alpha',\beta'}:(\lambda,v)\mapsto(\lambda,u)$,  is a
bijection from $\mathcal{F}_n(\alpha,*,\beta)$ to $\mathcal{F}_n(\alpha',*,\beta')$.

Combining the bijection  $\tau_{\alpha,\beta,\alpha',\beta'}$  with the involution $\Phi$ we
can prove bijectively the super-symmetry of  the polynomial $\sum_{(\lambda,v)\in \mathcal{F}_n} x^{a_v} y^{l_v}$.
More precisely we have the following result.
\begin{thm}
If $\alpha+\beta=\alpha'+\beta'$, the mapping $\zeta_{\alpha,\beta,\alpha',\beta'}=\Phi\circ\tau_{\alpha,\beta,\alpha',\beta'}\circ\Phi$ is a bijection from $\mathcal{F}_n(\alpha,\beta,*)$ to $\mathcal{F}_n(\alpha',\beta',*)$.

\end{thm}
\begin{proof}
Fix nonnegative integers $\alpha$, $\beta$, $\alpha'$, and $\beta'$ satisfying $\alpha+\beta=\alpha'+\beta'$.
By Theorem~\ref{thm:main},
the mapping $\Phi$ is a bijection from $\mathcal{F}(\alpha, \beta, *)$ to $\mathcal{F}_n(\alpha,*,\beta)$
and also a bijection from $\mathcal{F}(\alpha', *, \beta')$ to $\mathcal{F}_n(\alpha', \beta', *)$.
Since the mapping $\tau_{\alpha,\beta,\alpha',\beta'}$ is a bijection from $\mathcal{F}(\alpha, *, \beta)$ to $\mathcal{F}_n(\alpha',*,\beta')$,
it is obvious that the mapping $\zeta_{\alpha,\beta,\alpha',\beta'}=\Phi\circ \tau_{\alpha,\beta,\alpha',\beta'} \circ \Phi$ is a bijection from $\mathcal{F}_n(\alpha,\beta,*)$ to $\mathcal{F}(\alpha', \beta', *)$.
The bijection $\zeta_{\alpha,\beta,\alpha',\beta'}$ is illustrated as follows:
$$
\vcenter{
\xymatrix{
\mathcal{F}_n(\alpha,\beta,*) \ar[rr]^-{\zeta_{\alpha,\beta,\alpha',\beta'}} \ar[d]_{\Phi}
&&\mathcal{F}_n(\alpha',\beta',*)\\
\mathcal{F}_n(\alpha,*,\beta) \ar[rr]^-{\tau_{\alpha,\beta,\alpha',\beta'}} &&\mathcal{F}_n(\alpha',*,\beta') \ar[u]_{\Phi}
}
}.
$$
We are done.
\end{proof}
Theorem~4  yields that the generating function of $\mathcal{F}_n$  by the bivariate joint distribution of arm length and leg length is super-symmetric.
\section*{Acknowledgement}
This work is supported by la R\'egion Rh\^one-Alpes through the program ``MIRA Recherche 2008'', project 08 034147 01.

\bibliographystyle{amsalpha}

\providecommand{\bysame}{\leavevmode\hbox to3em{\hrulefill}\thinspace}
\providecommand{\MR}{\relax\ifhmode\unskip\space\fi MR }
\providecommand{\MRhref}[2]{%
  \href{http://www.ams.org/mathscinet-getitem?mr=#1}{#2}
}
\providecommand{\href}[2]{#2}

\end{document}